\documentclass[12pt]{amsart}
\usepackage{latexsym,fancyhdr,amssymb,color,amsmath,amsthm,graphicx,lettrine,listings,comment,anyfontsize,fancyhdr,tikz}
\newtheorem*{satz}{Theorem}
\definecolor{brightblue}{rgb}{0.9,0.9,1}
\definecolor{brightgreen}{rgb}{0.8,1.0,0.8}
\definecolor{brightyellow}{rgb}{1,1,0.7}
\def\summary#1{ \begin{center} \fcolorbox{black}{brightblue}{ \vspace{1mm} \parbox{15.2cm}{#1}} \end{center} }
\def\chapter#1{ \fcolorbox{brightgreen}{brightgreen}{ \parbox{15.2cm}{{\Large \begin{center}{\bf #1} \end{center} }}} }
\def\satz#1{ \vspace{1mm} \begin{center} \fcolorbox{black}{brightyellow}{ \parbox{15.2cm}{{\bf Theorem:} #1}} \vspace{1mm} \end{center} }
\let\paragraph\subsection

\makeindex

\title{Elements of Finite Geometry I}
\author{Oliver Knill}
\date{7/31/2026}
\address{Oliver Knill, knill@math.harvard.edu, Elements of Finite Geometry I, 2026}

\begin{document}
\thispagestyle{empty}
\begin{center}

    {\Huge\bfseries ELEMENTS OF } \vspace{1em} \\
    {\Huge\bfseries FINITE GEOMETRY I} \vspace{1em} \\

    \vspace{2em} \rule{0.1\linewidth}{0.5pt} \vspace{1em}

    {\large OLIVER KNILL}

    \vspace{0.5em}
    \vspace{8em}

\begin{center}
\scalebox{1.5}{\includegraphics{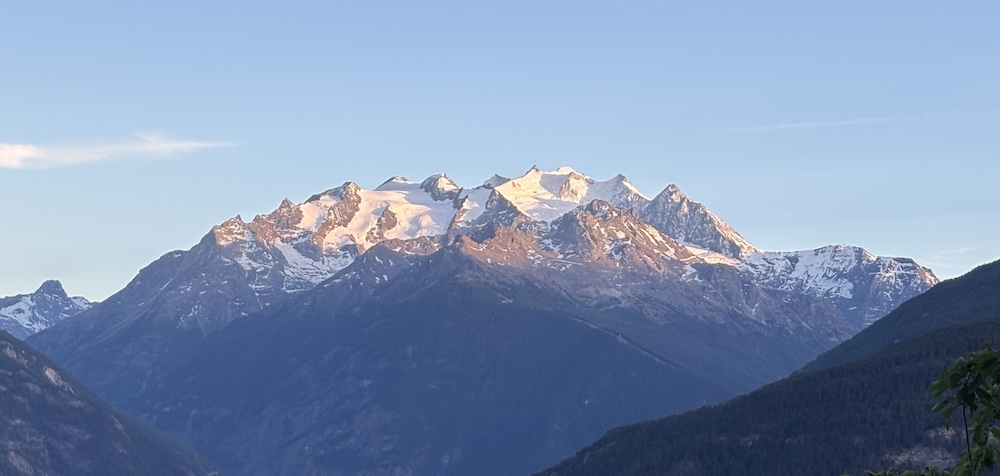}}
\end{center}

    \vspace{4em}
    \vfill
 \begin{tikzpicture}[scale=0.7]
        \fill[red] (0,0)--(4,0)--(2,3.4)--cycle;
        \fill[yellow] (0,-0.6) rectangle (4,0.0);
        \fill[blue]  (0.5,-0.2) arc[start angle=0,end angle=360,radius=0.5];
    \end{tikzpicture}

    \vspace{1em}

    {\large 2026}

\end{center}

\pagebreak

\maketitle

\setcounter{section}{0}  \chapter{Introduction}

\hspace{4mm}

\lettrine[lines=3]{T}{}here is a fundamental possibility that the common axiom system in mathematics
is inconsistent. A strong enough foundation of mathematics like the Zermelo-Fraenkel frame work can never be 
proven to be consistent within itself and is always extendable in the sense that there are statements that can 
not be proven nor dis-proven within the system. Hilbert's dream of building
a provably consistent strong theory can not be realized. The incompleteness theorems have entered popular culture 
especially with Hofstadter's book ``G\"odel-Escher-Bach". While most mathematicians - including myself - do not worry
about a possible inconsistency dooms-day, the possibility of such an event is a strong motivation to pursue 
geometry using weaker tools. Working with finite mathematics is one of these approaches. Finite descriptions of 
spaces like manifolds have already been initiated in Poincar\'e's ``analysis situs" program, but most topologists still
assume objects realized in an Euclidean space.
Even Dehn and Heegaard who initiated abstract simplicial complexes in 1907 used the continuum. 
\index{Zermelo-Fraenkel}
\index{Dehn-Heegaard}
\index{Hilbert's dream}

\hspace{2mm}

\lettrine[lines=1]{F}{}inite mathematics does not necessarily escape G\"odel. But it is conceivable 
that ZF is inconsistent while an axiom system without infinity will remain consistent. 
When working in a finite geometric set, we can in principle be certain to have consistency within that world. 
Most computer scientist are by nature more radical. They are ``strict finitists" in the sense that they only work with
structures for which there are tools to handle them. We can make sense of the symbol $GP = 10^{(10^{100})}$, 
but not of the number ``googolplex" itself. As computer scientists, we am essentially blind for answering questions like 
how many primes there are in $[GP-100,GP+100]$. Working with finite structures puts mathematics on safer grounds. It is a
refuge to which we might have to retreat one day, if some inconsistency in ZF
should emerge. I personally think this can happen but that there is little to worry: finite tools
will be able to emulate all continuum math. It is still possible that even after removing the infinity axiom, 
a catastrophic collapse could occur, but that is less likely; it would mean that even the algebra of finite sets 
is inconsistent. Assuming consistency of finite systems requires some faith, similarly as we trust the
Turing-Church thesis or assume that our memory in the brain is reliable, allowing us to do logical steps
without forgetting previous assumptions, while doing logical conclusions. In the last couple of years we all have seen 
intrinsically random probabilistic thinking entities to emerge that have started to compete with human thought. Aristotle
would turn in his grave learning that humanity started to flip coins in order to draw conclusions. Probabilistic reasoning is good 
enough for practical purposes, but it is philosophically unacceptable, even if the probability to be wrong is only $1/GP$. 
\index{finitism}
\index{googolplex}
\index{strict finitism}

\vspace{2mm}

\lettrine[lines=1]{A}{} ``finite geometry" is a finite geometric structure like a simplicial complex, 
a graph or a delta sets. These combinatorial structures carry natural
arithmetic and topological operations. One can add or multiply such structures in various ways for example.
The choice of using these three structures a bit arbitrary but it is motivated by the Unix trilogy 
"simplicity, clarity and generality": simplicial complexes have only one axiom and so are simple, graphs are 
intuitive and visual and so very clear. Delta sets finally are the most general. Every graph defines
a simplicial complex and every simplicial complex is naturally a delta set. Delta sets are sets
of sets that admit a calculus structure in the form of face maps. Their advantage of the delta set 
category is that has no limitations but still allows to do all geometry that we know from in the continuum. 
For example: submanifolds or products of geometric structures are often a
delta set at first before modeled again as a graph. Multi-graphs like quivers are delta sets. We work in this trinity. 
\index{delta set}
\index{graph}
\index{simplicial complex}

\vspace{2mm}

\lettrine[lines=1]{I}{}t is tempting to label a discussion about finite geometry with terms like ``quantum geometry".
The meaning of the word ``quantum" has been washed out recently and so has become ambiguous, 
even misleading.  The term ``quantum gravity" illustrates it. Classical geometry or calculus can be ``quantized" 
in various ways: there are non-commutative geometry flavors, Poisson brackets can be replaced by commutators or
Lagrange variational problem are deformed to path integral quantization. 
Every numerical scheme is some sort of quantization, as the continuum is modeled using finite structures 
suitable for numerical purposes. 
\index{quantum gravity}
\index{quantum geometry}
\index{quantum calculus}

\vspace{2mm}

\lettrine[lines=1]{F}{}inite geometry is related to ``calculus without limits". For example, 
every function $f$ on the real line satisfies the exact Taylor formula $f(a+x) = \sum_{n=0}^{\infty} f^{(n)}(a) x^n/n!$,
for integer $x$, if the derivative is deformed to $Df(a) = f'(a)=f(a+1)-f(a)$ and the polynomial algebra is deformed to 
$x^n = x (x-1) \dots (x-n-1)$. No regularity whatsoever is needed and Taylor is a finite sum. For $x=2$ for example 
$f(2)=f(0) + f'(0) 2/1! + f''(0) 2 (2-1)/2! = f(0) + 2f(1)-2f(0) + (f(2)-2 f(1) + f(0)) = f(2)$. If 
the unit $1$ is the Planck scale, finite calculus produces the same physics as traditional calculus. 
The derivative $D$ generates translation because $\exp(D t)f(x)=f(x+t)$ and $X f(x)=x f(x-1)$ and leads with 
the momentum operator $P=iD$ to the canonical anti-commutation relations $[X,P] = (XP-PX)=i$. 
\index{calculus without limits}
\index{Planck scale}
\index{Taylor formula}
\index{momentum operator}
\index{anti-commutation relations}

\vspace{2mm}

\lettrine[lines=1]{T}{}heorems are at the heart of mathematics. A theory needs to be judged
by the quality of theorems it can produce. Ideally, definitions are short and clear, theorems
are simple, proofs are intelligible and examples are rich. In a first batch we cover 12 results in finite geometry. 
I hope to be able to continue with the next batch of 12 at a later time.

 \vfill \pagebreak
\setcounter{section}{1}  \chapter{Unit 1: Gauss-Bonnet}

\summary{
The Gauss Bonnet theorem equates the total curvature of a geometry with its
Euler characteristic. Curvature is the energy of the fundamental units of 
space pushed to points.  }

\paragraph{}
A simple model for geometry is a {\bf finite simple graph} $(V,E)$. Its complete subgraphs 
define a finite abstract simplicial complex $G$,
a finite set of non-empty sets that is closed under the operation of taking non-empty subsets. 
The elements in $G$ are known as {\bf simplices}, {\bf cliques} or {\bf faces}.
The empty set $\emptyset$ is not a simplex. But the empty set is called the {\bf void} and is a
simplicial complex. A simplex $x$ of $(n+1)$ points has {\bf dimension} ${\rm dim}(x)=n$. 
Its {\bf energy} is $\omega(x) = (-1)^{{\rm dim}(x)}$. The total energy 
$\chi(G) = \sum_{x \in G} \omega(x)$ is the {\bf Euler characteristic} $(V,E)$. 
\index{Finite simple graph}
\index{Graph}
\index{simplicial complex}
\index{energy}
\index{Euler characteristic}
\index{void}

\paragraph{}
Let $f_k(G)$ denote the number of elements $x \in G$ with dimension $k$. For example, 
$f_0(G)=|V|$ is the number of vertices $K_1$, $f_1(G)=|E|$ the number of edges $K_2$
and $f_2(G)$ is the number of triangles $K_3$ in the graph. 
We have $\chi(G) = \sum_{k=0}^q (-1)^k f_k(G)$, where $q$ is the {\bf maximal dimension.}
The vector $(f_0(G),f_1(G), \dots, f_q(G))$ is called the $f$-vector of $G$. The number 
$q$ is the {\bf maximal dimension} of the graph.
\index{f-vector}
\index{maximal dimension}

\paragraph{}
The {\bf unit sphere} of a vertex $v$ is the sub-graph generated by all its neighbors:
$S(x) = (W=\{ w \in V, (v,w) \in E \}, \{ e=(a,b) \in E, a \in W, b \in W \} \})$. It
has maximal dimension $\leq q-1$ if $G$ has maximal dimension $q$.
The {\bf curvature of a vertex} is $v$ is defined as 
$$  K(v) = 1-\sum_{k=0}^{q-1} \frac{ (-1)^{k} f_k(S(v))}{k+1} \; . $$
\index{Curvature}
\index{Unit sphere} 

\satz{
$\chi(G)=\sum_{v \in V} K(v)$. 
}
\begin{proof}
Distribute the energy $w(x)$ of each simplex $x$ of dimension $k$ equally 
to each of its $k+1$ vertices. Then $K(v)$ is the total energy, that $v$ has collected. 
\end{proof} 
\index{Gauss-Bonnet theorem}

\paragraph{}
A different proof is obtained by noting the $\chi(G)$ is a linear combination
of valuations $f_k(G)$ which satisfy $f_k(G) = \frac{1}{k+1} \sum_{v} f_{k-1}(S(v))$
generalizing the {\bf Euler handshake formula} $f_1(G) = \frac{1}{2} \sum_v f_0(S(v)) = \frac{1}{2} \sum_v d(v)$,
where $d(v)=f_0(S(v))$ is the {\bf vertex degree} of $v$.
\index{Handshake formula}
\index{valuations}
\index{vertex degree}

\paragraph{}
To generalize the theorem to simplicial complexes $G$, let $V$ be the set of $0$-dimensional simplices in $G$.
$V$ is naturally identified with $\bigcup_{x \in G} x$. Given $x \in G$, define $U(x)=\{ y \in x, x \subset y\}$.
The unit sphere is $S(v) = \{ x, v \in x, \{v \neq x\} \}= \overline{U(v)}\setminus  U(v)$
is a simplicial complex for which $f_k(S(v))$ is defined. Curvature is supported on $V$.

\paragraph{}
The {\bf simplex generating function} $f_G(t)=1+f_0(G) t + \cdots + f_d(G) t^{q+1}$ 
of $G$ satisfies $\chi(G)=1-f_G(-1)$.
The derivative of a polynomial can be defined without taking limits by
using the rule $f'=\sum_{n=0}^q a_n n x^{n-1}$ if $f=\sum_{n=0}^q a_n x^n$.
With the {\bf curvature generating function} $K_G(t) = \int_0^t f_G(s) \; ds$, 
we have $K(v)=K_{S(v)}(-1)$. Gauss-Bonnet is now an elegant recursive formula:
\index{index generating function}
\index{curvature generating function}

\satz{ $f_G'(t) = \sum_{v \in V} f_{S(v)}(t)$.  }

\begin{proof}
Use the Gauss-Bonnet theorem and check that
$K(v) = \int_{-1}^0 f_{S(v)}(t) \; dt$. 
\end{proof}
\index{functional Gauss-Bonnet}

{\bf Examples}:  \\

1) If $G$ is $1$-dimensional, there are no triangles 
and $\chi(G)=f_0(G)-f_1(G)=|V|-|E|$. The curvature is $K(v) = 1-{\rm deg}(v)/2$,
where ${\rm deg}(v)$ is the {\bf vertex degree}.  \\
2) {\bf Circular graphs} $C_n$ for example have curvature constant $0$.  \\
3) {\bf Star graphs} $S_n$ have curvature $1/2$ on the $n$ leafs and curvature $1-n/2$ 
at the branch center. The total curvature of a star graph is $1$.  \\
4) {\bf Path graphs} $P_n$ have curvature $1/2$ at the boundary and $0$ else. 
The {\bf cube graph} $C$ has constant curvature $K(v) = 1-3/2=-1/2$ leading to $\chi(C)=-4$.
5) The {\bf dodecahedron graph} $D$ has $K(v)=1-3/2=-1/2$ and $\chi(D)=-10$.  \\
\index{Octahedron}
6) A graph is called a {\bf $2$-manifold}  if every 
unit sphere $S(v)$ is a circular graph of length $\geq 4$. The Euler characteristic is then 
$\chi(G) = |V|-|E|+|F|$. An example is the octahedron $O$ with $|V|=6$ vertices
or the icosahedron $I$ with $|V|=12$ vertices. The curvature of a 2-dimensional manifold is 
the {\bf Eberhard formula} $1-{\rm deg}(v)/6$ because
$f_0(S(v)) = f_1(S(v)) = {\rm deg}(v)$ and $f_0(S(v))/2 - f_1(S(v))/3={\rm deg}(v)/6$.  
For the octahedron, the curvature is constant $K(v)=1-4/6=1/3$ adding up to $\chi(G)=2$.
For the icosahedron, the curvature is constant $K(v)=1-5/6=1/6$ adding up to $\chi(G)=2$.  \\
7) A complete graph $K_n$ has $n$ vertices and constant $K(v)=1/n$. 
One can see this by noting that the curvature must be constant and add up to $1$
and that the Euler characteristic of $K_n$ is $1$. Note that 
$f_k(K_n)=\left( \begin{array}{c} n \\ k \end{array} \right)$
so that $(1+t)^n = \sum_{k=0}^n f_{k-1}(K_n) t^k$ and
for $t=-1$, this is $0$ to the left and $-1+\chi(K_n)$ on the right. 
For the {\bf tetrahedron} $T=K_4$ for example, the curvature is constant $1/4$. It is not a
2-manifold.   \\
8) The 4-partite graph $K_{2,2,2,2}$ is the {\bf 3-dimensional cross polytope}. It is a
discrete {\bf 3-sphere} or $16$-cell. Its $f$-vector is $f=(f_0,f_1,f_2,f_3)$ $=(|V|,|E|,|F|,|C|)$ 
$=(8,24,32,16)$, where $F$ is the set of triangles, the $2$-dimensional parts of space and 
$C$ is the set of tetrahedra, the $3$-dimensional parts of space. The Euler characteristic is
$\chi(G)=8-24+32-16=0$, like for any odd-dimensional manifold. The curvature is constant $0$.

\index{triangle free graphs}
\index{Circular graphs}
\index{Star graphs}
\index{Path graph}
\index{cube graph}
\index{Eberhard formula}
\index{2-manifold}
\index{Complete graph}
\index{Tetrahedron}
\index{Cross polytope}
 \vfill \pagebreak
\begin{center}\scalebox{0.8}{\includegraphics{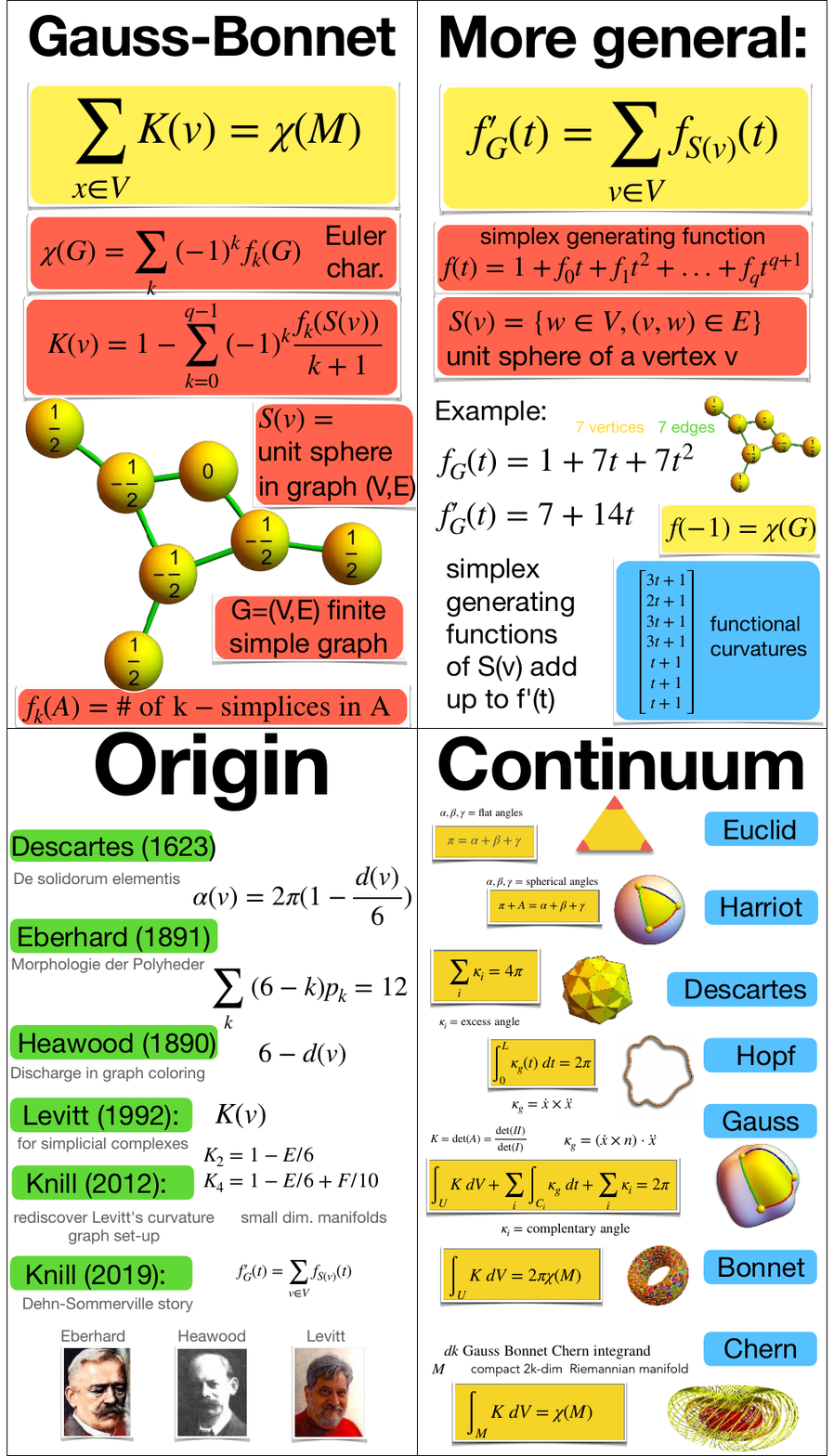}}\end{center} \pagebreak
\setcounter{section}{2}  \chapter{Unit 2: Poincare-Hopf}

\summary{
The Poincare-Hopf theorem writes the Euler characteristic of
a space as the degree of a divisor, the sum of integer indices attached to points.
This works for locally injective scalar function or locally non-circular vector fields.  }

\paragraph{}
Geometry is again a simplicial complex $G$ from a finite simple graph $(V,E)$.
A {\bf scalar function} $g$ is a map $g:V \to R$, where $R$ is a totally ordered
set, like a subset of $\mathbb{Z},\mathbb{Q}$ or $\mathbb{R}$. We assume that $g$ is a {\bf coloring}, 
meaning that it is {\bf locally injective}:
$g(v) \neq g(w)$ if $(v,w) \in E$. For $v \in V$,
the {\bf stable sphere} $S^-_g(v)$ is the sub-graph of the unit sphere $S(v)$,
generated by $\{ w \in S(v), g(w)<g(v) \}$. Define the {\bf index} $i_g(v) = 1-\chi(S^-_g(v))$. 
\index{stable sphere}
\index{coloring}
\index{locally injective}
\index{index}

\satz{ $\chi(G)=\sum_{v \in V} i_g(v)$.  }
\index{Poincar\'e-Hopf}
\begin{proof}
Distribute the energy $\omega(x)=(-1)^{{\rm dim}(x)}$ of each simplex $x \in G$ to the vertex $w$ in $x$, 
where $g$ takes its maximum. Then $i_g(v)$ is the total energy which the 
vertex $v$ has collected because every $k$-simplex $x$ in $S^-_g(v)$ gives its energy $\omega(x)$ to $v$.
\end{proof} 

\paragraph{}
For the {\bf simplex generating function} $f_G(t)=1+f_0(G) t + \cdots + f_d(G) t^{d+1}$, this means:
\index{simplex generating function}

\satz{
$f_G(t) = 1+t \sum_{v \in V} f_{S^-_g(v)}(t)$.
}
\begin{proof}
The k-simplex $x$ in $S(v)$ matches up with a $(k+1)$-simplex in $G$ that
contains $v$. This shift of dimension is the reason for the multiplication with $t$. 
\end{proof} 

\paragraph{}
This formula again gives a recursive computation of $f_G(t)$ and so the cliques of $G$. 
It is even more efficient than Gauss-Bonnet because $S^-_g(v)$ is in general smaller 
than $S(v)$. The index of $v$ is a special value of the generating function of $S_g^{-1}(v)$:
$$ i_g(v) = f_{S_g^-(v)}(-1)  \; . $$

\paragraph{}
A {\bf vector field} is a map $F: G \to V$ with $F(x) \in x$. 
A {\bf directed graph} with the property that there are no closed loops in 
each simplex $x$ defines a vector field by $F(x)={\rm max}_{v \in x}(v)$,
as the directions on each edge define a total order on $x$. We call 
this a {\bf locally non-circular digraph}. 
\index{vector field}
\index{directed graph}

\paragraph{}
Define $S^-_F(v) = \{ x \in G, F(x)=v \} = F^{-1}(v)$ and
$i(v) = 1-\chi(S_F^-(v))$. We again have
$f_G(t) = 1+t \sum_{v \in V} f_{S^-_F(v)}(t)$.

\paragraph{}
An example is the
{\bf gradient field} defined by a locally injective function $g$.
In this case $F(x)=v$, where $v$ is the maximum of $g$ on $x$. 
The index of such a graph is then 
$i(v) = 1-\chi(S^-(v))$, where $S^-(v)=\{ w \in S(v), w<v \}$. 
\index{gradient field}

\satz{ For a locally non-circular directed graph, $\chi(G)=\sum_{v \in V} i(v)$.  }
\index{Poincar\'e-Hopf for directed graphs}

\paragraph{}
Define the {\bf symmetric index} $j_g(v) = [i_g(v)+i_{-g}(v)]/2$. By linearity,
we still have $\sum_{v \in V} j_g(v)=\chi(G)$. 
\index{symmetric index}

\begin{enumerate}
\item Let $G=C_n$ be a cycle graph and $g$ a function. Now
$i_g(v)=1$ for local minima and $i_g(v)=-1$ for local maxima. Since 
minima and maxima alternate, there are the same number of minima and maxima. 
\item If $G=S_n$ is a star graph, and the maximum is at the center $c$ then 
$i_g(c)=n-1$ and all leaves have $i_g(v)=1$. The total of all index values is $1$. 
If $g(c)$ is the $k$'th maximal entry then $i_g(c)=n-k$ and $k-1$ entries
that are smaller produce leaf indices $1$. 
\item  If $v$ is a local minimum of $g$ on $V$, then $i_g(v)=1$ because
$S(v)=0$ is the empty graph which has Euler characteristic $0$.
If $v$ is a local maximum of $g$ on $V$, then $i_g(v)=1-\chi(S(v))$.
For a 2-manifold for example where $S(v)$ is a circular graph we have 
$i_g(v)=1$ at a maximum too. At a saddle point, where $S(v)$ consists of 
two disconnected arcs, we have $i_g(v)=-1$. 
\item If $G$ is a complete graph and $g$ is an arbitrary function, then 
$i_g(v)=1$ at the minimum and $i_g(v)=0$ else. The symmetric index is $1/2$ on 
the maximum and $1/2$ on the minimum. 
\item If $G$ is a manifold for which $\chi(S(v))=2$ for all $v$ like
for odd dimensional manifolds, then $\chi_{S^-_{-g}}(v)=2-S^-_{g}$ which can 
be written as $i_g(v)=-i_{-g}(v)$ so that $j_g(v)=0$ for all $v$. 
Odd dimensional manifolds all have Euler characteristic $0$, independent of whether
they are orientable. 
\item If $G_1$ is the Barycentric refinement of $G$ in which the vertices are
the simplices of $G$ and two vertices are connected if one is contained in the other, 
look at the coloring $g(x)={\rm dim}(x)$. Now $S^-_g(x)$ is the boundary complex of the simplex
$x$ and $1-\chi(S^-_g(x))=\omega(x)$. The Poincar\'e-Hopf formula now tells
that $\chi(G)=\sum_{x} \omega(x) = \sum_{v \in V(G_1)} i_g(v) = \chi(G_1)$. 
The Barycentric refinement $G_1$ has the same Euler characteristic than $G$. 
\end{enumerate}
 \vfill \pagebreak
\begin{center}\scalebox{0.8}{\includegraphics{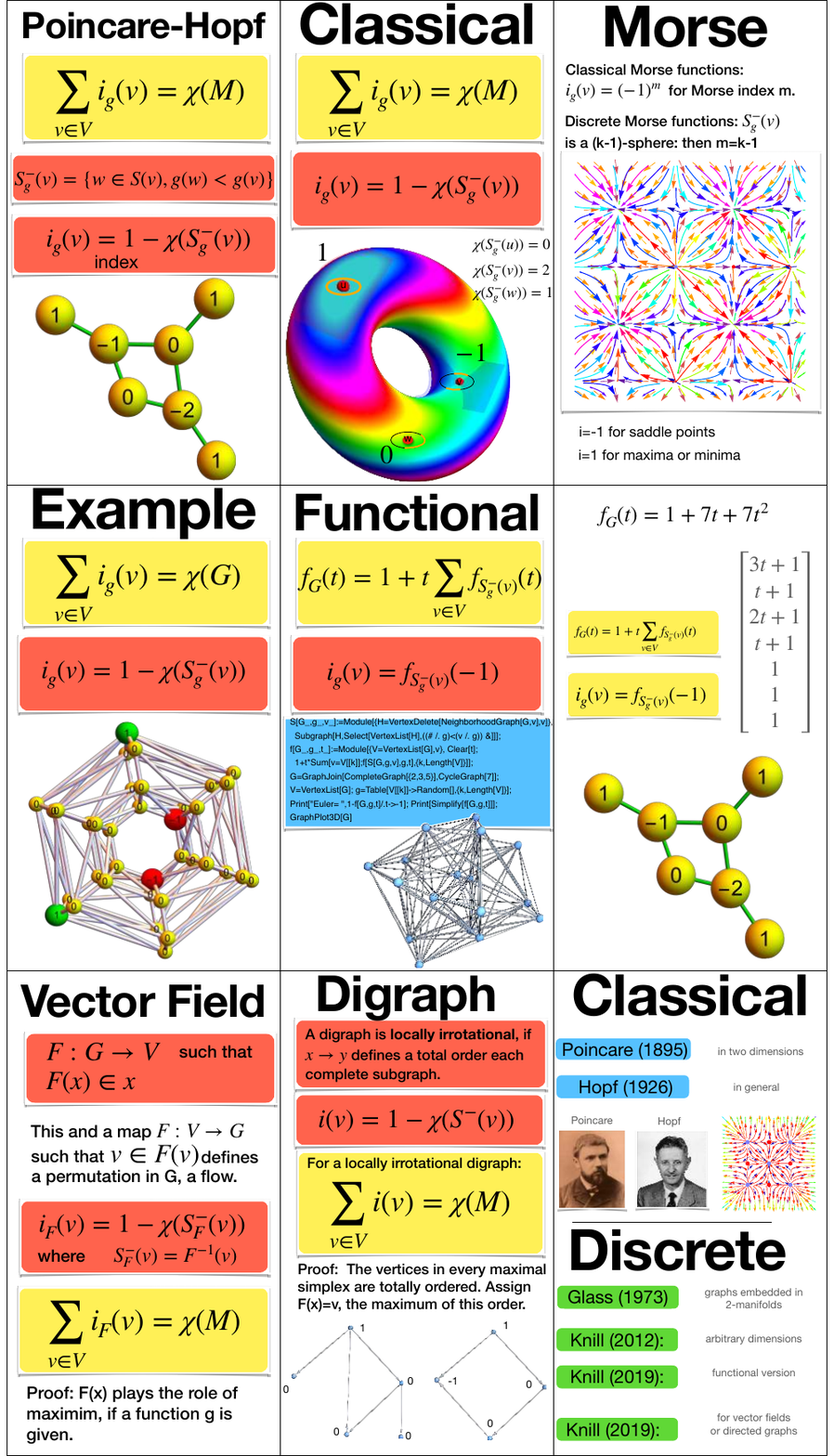}}\end{center} \pagebreak
\setcounter{section}{3}  \chapter{Unit 3: Index Expectation}

\summary{
Integral geometry links Poincar\'e-Hopf and Gauss-Bonnet: curvature is
the expectation of indices. This is a tool that can link the discrete with the continuum.  }

\paragraph{}
Let $G$ be a simplicial complex from a $(V,E)$ is a finite simple graph and let $R$ 
be a totally ordered set, like an interval $[a,b]$ or a subset of the integers. Consider 
$(\Omega,\mathcal{A},{\rm P})$, a {\bf probability space} of locally 
injective functions from $V$ to $R$. An example is
the set of all possible colorings from $V$ to $\{1,\dots, c\}$, where $c$ is 
the {\bf chromatic number} and assume that ${\rm P}$ is the counting
measure giving each coloring the same probability. For fixed $v$, 
the index $g \to i_g(v)$ is now an integer valued {\bf random variable} on $\Omega$. 
\index{Probability space}
\index{coloring}
\index{chromatic number}
\index{random variable}

\paragraph{}
The {\bf index expectation} $K(v)={\rm E}[i_g(v)]$ does not involve $g$ any more.
it defines a curvature that depends on the probability space.
\index{index expectation}
\index{curvature as index expectation}

\satz{
If $K(v)={\rm E}[i_g(v)]$, then 
$\sum_{v \in V} K(v) = \chi(G)$. 
}
\begin{proof}
Start with Poincar\'e-Hopf $\sum_{v \in V} i_g(v) = \chi(G)$. 
Take the expectation of this formula to get
${\rm E}[\sum_{v \in V} i_g(v)] = {\rm E}[\chi(G)] = \chi(G)$. 
Fubini (or rather linearity of expectation) allows to write the left hand side as 
$\sum_{v \in V} {\rm E}[i_g(v)] = \sum_{v \in V} K(v)$. 
\end{proof} 
\index{Fubini}

\paragraph{}
Every $v \in V$ defines a random variable $X_v(g) = g(v)$ on $(\Omega,\mathcal{A},{\rm P})$.
If these random variables are independent and have all the same uniform 
distribution on some interval, we say ${\rm P}$ has the {\bf Lebesgue property}.
\index{Lebesgue property}

\satz{If ${\rm P}$ has the Lebesgue property, then $K$ is the Levitt curvature.}
\begin{proof}
For every $k$-simplex $x=(x_0, \dots, x_k)$, the probability that a function 
$g$ has the maximum on $x_j$ is the same as the values of $g$ outside $g$ are 
independent and so do not factor in. By assumption, the situation on the simplex
is now symmetric with respect to any permutation of the vertices in simplices. This means
that each vertex gets $1/(k+1)$ of the energy $\omega(x)$.
\end{proof} 

\paragraph{}
Assume $c$ is the {\bf chromatic number} of $G$. Let $\Omega$ denote the set of 
all $c$ colorings. If ${\rm P}$ has the uniform distribution on this finite set, 
we say ${\rm P}$ is the {\bf coloring probability space}.

\satz{For the coloring probability space, $K$ is the Levitt curvature.}
\begin{proof}
The same symmetry argument applies: as the probability measure does not change
if we permute the coloring space, each point $x_j$ in a simplex has the same
probability of being the maximum.
\end{proof}

\paragraph{}
The index itself is  curvature if ${\rm P}$ is supported
on a single point $g$ in $\Omega$. 
By changing the ${\rm P}$, we can {\bf deform} a geometry. Index expectation
allows us to modify the curvature similarly as a Riemannian metric does in 
the continuum. 
\index{deforming a geometry}

\paragraph{}
If the graph $(V,E)$ is a triangulation of a compact Riemannian manifold $M$ 
that is isometrically Nash embedded into an ambient Euclidean space $E$, then 
almost all linear functions $g(x) = x \cdot a$ induce
Morse functions on $M$. For a sufficiently fine triangulation $G$ in $M$, almost 
all linear functions in the ambient space are locally injective on $G$. 
The Poincar\'e-Hopf indices at critical points match to 
the classical ones. The index expectation is a curvature that is locally 
homogeneous. An argument of Weyl implies that it 
has to be the Gauss-Bonnet-Chern integrand.
\index{Nash embedded}

\paragraph{}
Curvature is compatible with the {\bf Shannon product} $A*B$ of two graphs $A=(V,E),B=(W,F)$ is
$(V \times W, Q)$ with edges $Q=\{ ( (a,b),(c,d) ), {\rm dist}(b,d) \leq 1 \; {\rm and} \; {\rm dist}(a,c) \leq 1 \}$
where ${\rm dist}$ stands for the {\bf graph distance}. It is also called {\bf strong product}
\index{Shannon product}
\index{strong product}

\satz{The curvature $K(x,y)$ on $A*B$ is $K(x,y) = K(x) K(y)$.}

\begin{proof}
(i) First show $i_{G*H,g h}(x,y) = i_{G,g}(x) i_{H,h}(y)$.
Proof: For any graphs $G,H$, we have $(1-\chi(G)) (1-\chi(H)) = (1-\chi(G \oplus H))$,
where $G \oplus H$ is the graph join. 
The formula $(1-\chi(S_{G*H,g*h}(x)))= (1-\chi(S_{G,g}(x))) (1-\chi(S_{H,h}(x)))$
follows now from $S_{G*H,g*h}(x,y)$ being homotopic to the join of
$S_{G,g}(x)$ and $S_{G,h}(y)$ and that the Euler characteristic is a homotopy invariant.
We have $S_{G*H}^-(x,y) = B_G^-(x)*S_H^-(y) \cup S_G^-(x)*B_H^-(y)$,
which is homotopic to the join of $S_H^-$ and $S_G^-$. \\
(ii) Take the expectation of the relation in (i)
and make use of the fact that the random variables $X(g)=i_{G,g}(x)$ and
$Y(h)=i_{H,h}(y)$ are independent. The expectation of the product is the product of the
expectations.  Since $K(x) = {\rm E}[X] ,  K(y) = {\rm E}[Y]$, we have
$K(x,y) = {\rm E}[X Y]  = {\rm E}[X] {\rm E}[Y]$.
\end{proof}

\paragraph{}
{\bf Examples}: 
1) If $A=K_n$ and $B=K_m$, then $A*B=K_{nm}$. The curvature on $A$ is constant $1/n$
the curvature on $B$ is constant $1/m$. The curvature on $A*B$ is constant $1/(nm)$. \\
2) Assume $G$ is embedded in $\mathbb{R}^2$, meaning $P(v)=(x(v),y(v))$ is given.
Let $g(v)  = \cos(\theta) x(v) + \sin(\theta) y(v)$ and $\Omega=\{ \theta \in [0,2\pi) \}$. 
Index expectation defines a curvature that depends on the embedding.  \\
3) Since the symmetric index $[i_g + i_{-g}]/2$ is zero for an odd dimensional manifold
any probability measure that is invariant under the involution $g \to -g$ produces 
a curvature that is constant zero on an odd-dimensional manifold. 
 \vfill \pagebreak
\begin{center}\scalebox{0.8}{\includegraphics{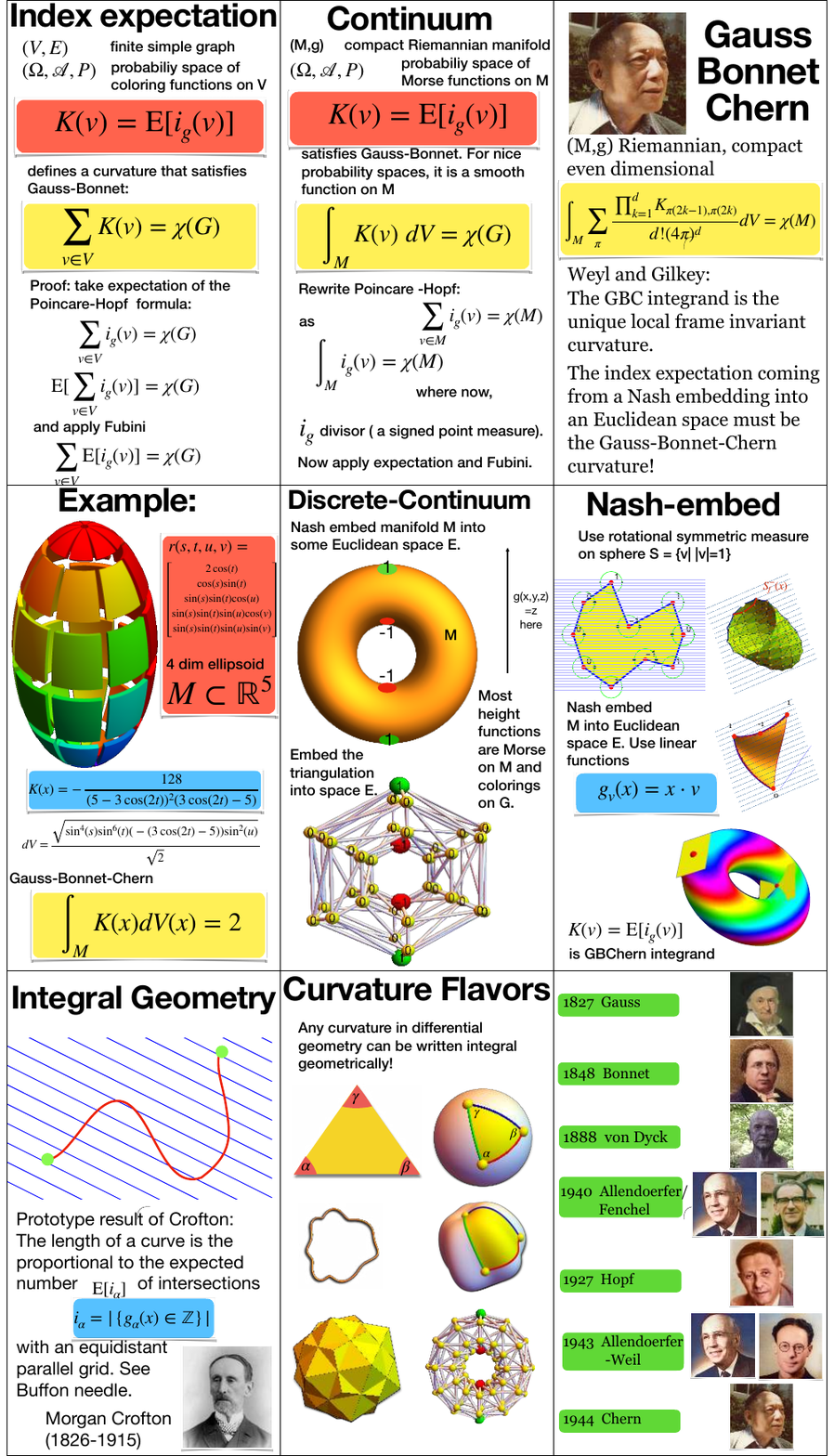}}\end{center} \pagebreak
\setcounter{section}{4}  \chapter{Unit 4: Euler's Gem}

\summary{
A $q$-sphere $G$ is a $q$-manifold which when punctured becomes
contractible. Euler's Gem formula tells that Euler characteristic is $\chi(G)=1+(-1)^q$.
Platonic q-spheres are inductively defined as $q$-spheres for which 
all unit spheres are isomorphic and equal to some Platonic $(q-1)$-sphere. They can be classified.
}

\paragraph{}
A graph $(V,E)$ is {\bf contractible} if there exists $v \in V$ such that
$S(v)$ and $G \setminus v$ are both contractible. The $1$-point graph $1=K_1$ is 
contractible. The $(-1)$-sphere $0$ is not contractible. 
\index{contractible}

\paragraph{}
A graph $G$ is called a {\bf q-manifold} if all unit spheres $S(v)$ are $(q-1)$-spheres.
A {\bf $q$-sphere} is a $q$-manifold for which there is a $v$ such that $G \setminus v$
is contractible. 
\index{q-manifold}
\index{q-sphere}

\paragraph{}
Every contractible graph has $\chi(G)=1$. One can see this by induction:
$\chi(G)=\chi(G \setminus v) + \chi(S(v)) - \chi(\{v\}) = 1+1-1=1$.
If $A,B$ are graphs, $A \cap B$ carries the intersection of the
simplicial complexes of $A$ and $B$, then the {\bf valuation formula} for 
simplicial complexes holds also for graphs $\chi(A) + \chi(B)-\chi(A \cap B) = \chi(A \cup B)$.
$A \cup B$ carries the intersection of the simplicial complexes of $A$ and $B$. 
The Euler Gem formula is:
\index{valuation formula}

\satz{
$\chi(G) = 1+ (-1)^q$ for a $q$-sphere.
}
\begin{proof}
(i) In every graph, every unit ball $B(v)$ is contractible. It follows that $\chi(B(v))=1$. 
(ii) for a $q$-sphere, every unit sphere $S(v)$ is a $(q-1)$-sphere and so
has by induction the Euler characteristic $1+(-1)^{q-1}) = 1-(-1)^q$. \\
(iii) Now use $\chi(G) = \chi(G \setminus v) + \chi(B(v)) -\chi(S(v))$ 
  $=1+1-(1+(-1)^{q-1}) = 1+(-1)^q$. We used the valuation formula.
\end{proof} 
\index{Euler's Gem}

\paragraph{}
The {\bf join} of two graphs $A \oplus B$ has vertex set $V(A) \cup V(B)$ 
and edge set $E(A) \cup E(B) \cup \{ (a,b), a \in V(A), b \in V(B) \}$. 
It is the dual of disjoint union $+$ because $A \oplus B = \overline{\overline{A} + \overline{B}}$. 
With this operation and the zero element $0$ we have a second monoid, dual and isomorphic
to the monoid with disjoint union. Spheres form a sub-monoid of the join monoid.
\index{join}
\index{join monoid}
\index{disjoint union monoid}

\satz{A is a $k$-sphere, $B$ is a $l$-sphere $\Rightarrow$ $A \oplus B$ is a $(k+l+1)$-sphere.  }
\begin{proof}
Use induction with respect to $k+l=n$. For $n=0$ meaning $k=l=0$, the join is a $1$-sphere. 
Assume we have proven it for all cases with $k+k \leq n-1$ and assume $k+l=n$. 
By definition every unit sphere $S_A(v)$ is a $(k-1)$-sphere and every unit sphere $S_B(w)$ is a $(l-1)$-sphere. 
For $v \in V(A)$, the graph $S_{A \oplus B}(v) = S_A(v) \oplus B$ is a $l-1+k+1$-sphere by induction.
For $w \in V(B)$, the graph $S_{A \oplus B}(w) = A \oplus S_B(w)$ is a $l+k-1+1$-sphere by induction.
This shows that $A \oplus B$ is a $l+k+1$ manifold.
Now use that if $A$ is contractible and $B$ is arbitrary than $A+B$ is contractible. 
This shows that for $v \in V(A)$, the graph $A \oplus B \setminus v = (A \setminus v)  \oplus B$ is contractible. 
Similarly, for $w \in V(B)$, the graph $A \oplus B \setminus w  = A \oplus (B \setminus w)$ is contractible. 
\end{proof}
\index{join of spheres}

\paragraph{}
Inductively, a q-sphere is declared to be {\bf Platonic}, if all unit spheres are Platonic $(q-1)$- spheres 
that are all isomorphic. $0$ is Platonic. There is a unique Platonic sphere in all dimensions except dimensions 1,2,3. 
There are infinitely many 1-spheres $C_n$ and they are all Platonic. 
For $d=2$,there is the {\bf octahedron} and {\bf icosahedron}.
For $d=3$, there is the {\bf 600 cell} and the {\bf 16 cell}.
After that there are only {\bf cross polytopes}, the $2^{q+1}$ cell. 
\index{Octahedron}
\index{Icosahedron}
\index{600 cell}
\index{16 cell}
\index{Platonic sphere}
\index{cross polytope}

\satz{ There is a unique Platonic sphere for $q>3$. }

\begin{proof}
$q=-1,0,1$ are clear. For $q=2$, the curvature 
$K(x)=1-f_0/2 + f_1/3-f_2/3$ is constant, adding up to $2$. 
It is either $1/3$ or $1/6$. For $d=3$, where
each $S(x)$ must be either the octahedron or icosahedron, 
$G$ is the 16 cell or 600 cell. For $d=4$, by Gauss-Bonnet, 
$K(x)$ add up to $2$ and be of the form $L/12$ for some integer $L$.
For $L=1$, there exists the $4$-dimensional cross polytope with 
$f$-vector $(10, 40, 80, 80, 32)$.
There is no $4$-sphere, for which $S(x)$ is the 600-cell. 
As the $f$-vector of it is $(120,720,1200,600)$, 
we would get $K(x)=1-120/2+720/3-1200/4+600/5 =1$. 
Gauss-Bonnet would give $|V|=2$ and ${\rm dim}(G) \leq 1$. 
Inductively there is now only one Platonic $q$-sphere for $q>3$. 
\end{proof} 

\paragraph{}
In classifications of {\bf Platonic solids}, {\bf dual polytopes} are usually included. Since they are
triangle free, they are not q-manifolds.
Classics also add the (q+1)-simplices because their boundary simplicial 
complex are $(q-1)$-spheres as simplicial complexes. The 24-cell $U$, the units of the 
Hurwitz quaternions is the only additional polytop. One could prove this
by extending Gauss-Bonnet to combinatorial cell-complexes: inductively add $k$-cells
to an already existing $(k-1)$ sphere. The notion of q-sphere, q-manifold, 
contractibility and Platonic sphere can be generalized to such
discrete complexes similarly as CW complexes in the continuum. After extending
Gauss-Bonnet in the same way, we build up the Platonic spheres inductively by dimension.
$U$ does not occur as a unit sphere of any Platonic $4$-complex so that the miracle in 
dimension $q=3$ ends. Traditionally, the continuum is involved:
regular q-polytopes are boundaries of convex regions in $\mathbb{R}^{q+1}$. 
\index{Platonic solids}
\index{dual polytopes}
\index{24 cell}

 \vfill \pagebreak
\begin{center}\scalebox{0.8}{\includegraphics{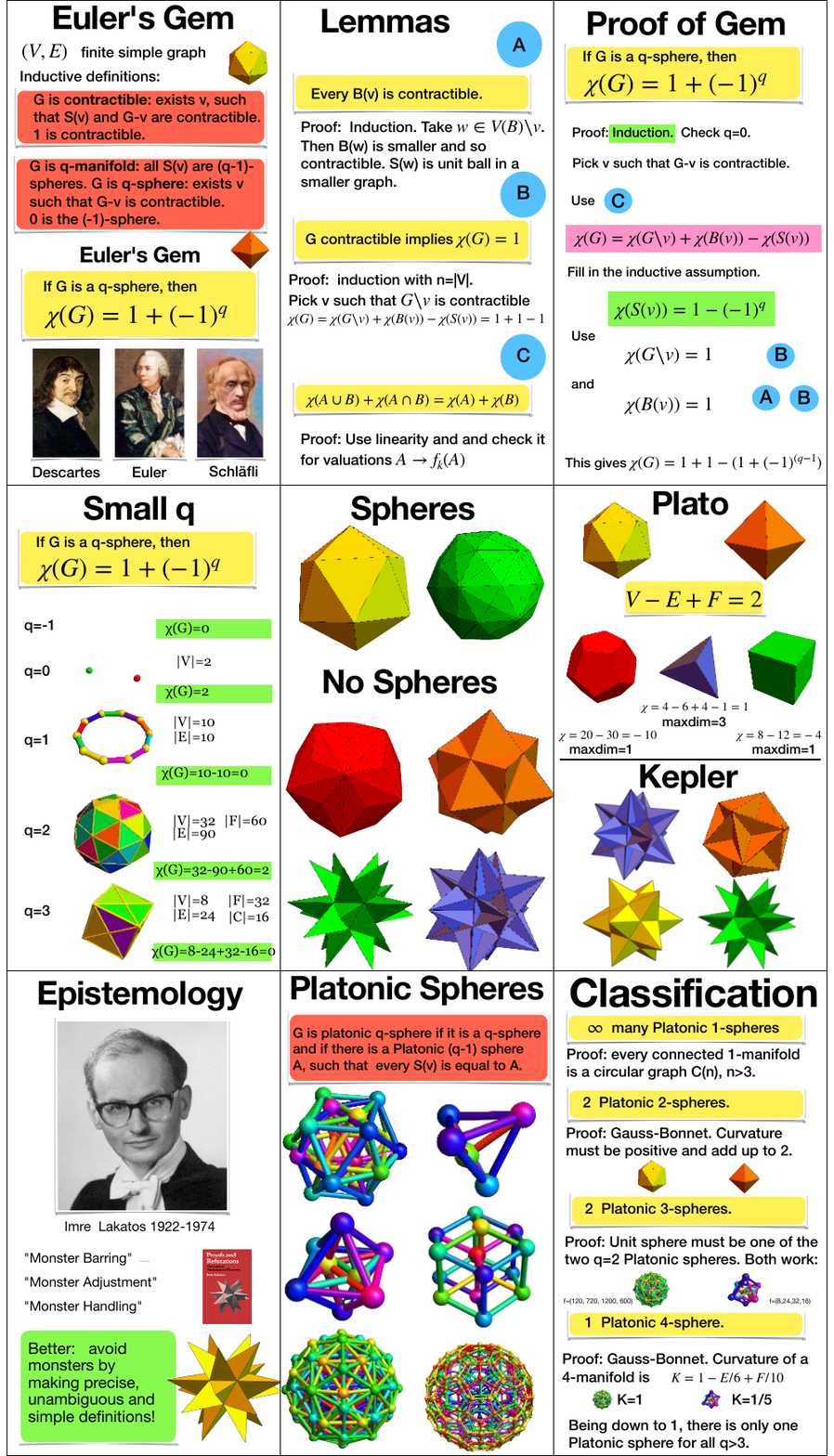}}\end{center} \pagebreak
\setcounter{section}{5}  \chapter{Unit 5: Euler-Poincar\'e}

\summary{
A simplicial complex defines an exterior derivative $d$
on the Hilbert space of forms. The kernels of $H=(d+d^*)^2$ restricted to $k$-forms is
the $k$'th Betti number. The alternating sum of the $b_k$ is the Euler characteristic.  }

\paragraph{}
Let $G$ be a finite abstract simplicial complex. It can be the Whitney complex of a graph for example.
If $G$ has cardinality $n$, the vectors in the Hilbert space 
$\Lambda=l^2(G) \sim \mathbb{R}^n$ are {\bf discrete differential forms} or simply {\bf forms}. 
Functions on $G_k = \{ x \in G, {\rm dim}(x)=k \}$ are called {\bf $k$-forms}. 
There is a natural decomposition $\Lambda = \oplus_{k=0}^q \Lambda_k$ and linear maps on $\Lambda$
that preserve $\Lambda_k$ are block matrices if we assume the elements in $G$ to be ordered according
to dimension.
\index{differential forms}
\index{forms}

\paragraph{}
For $y \subset x$ with ${\rm dim}(y)={\rm dim}(x)-1$, define $d(x,y) = {\rm sign}(x,y)$,
where ${\rm sign}(x,y)$ is the signature of the permutation that maps 
$(vy)$ to $x$, where $v=x \setminus y$. For example, 
${\rm sign}( (1,2,3,4),(1,2,4) ) = {\rm sign}( (1,2,3,4),(3,1,2,4) ) = 1$. The matrix $d$ is 
nilpotent $d^2=0$ and is called the {\bf exterior derivative}. It maps $k$-forms to $(k+1)$-forms. 
The symmetric matrix $D=d+d^*$ is the {\bf Dirac operator}. 
Its square $H=D^2 = d d^* + d^* d$ is the {\bf Hodge operator}. 
\index{signature}
\index{exterior derivative}
\index{Dirac operator}
\index{Hodge operator}

\paragraph{}
$H$ decomposes into blocks $\oplus_{k=0}^q H_k$ because both $d d^*$ and $d^* d$ map $\Lambda_k$
into itself. The {\bf $k$'th cohomology group} is defined as the space of {\bf $k$-harmonic forms} 
${\rm ker}(H_k)$. Its dimension is the $k$'th {\bf Betti number} $b_k$ of $G$. 
If $I_k$ are the identity matrices on $l^2(G_k)$, then 
$f_k={\rm tr}(I_k) = |G_k|$ is the number of $k$-simplices and $\chi(G)=\sum_{k=0}^{q} (-1)^k f_k$.  
Similarly as the f-vector $(f_0,f_1, \dots, f_q)$ encodes cardinalities, 
the {\bf Betti-vector} $(b_0,b_1 \dots, b_q)$ encodes dimensions of kernels. 
\index{harmonic forms}
\index{Betti number}
\index{Betti vector}

\paragraph{}
The {\bf McKean-Singer symmetry} is a spectral symmetry telling that non-zero eigenvalues on
even forms correspond bijectively to non-zero eigenvalues on odd forms. It can be encoded as:
\index{McKean-Singer symmetry}

\satz{ ${\rm str}(e^{-t H}) = \chi(G)$.  }
\begin{proof}
If $u$ is an eigenvector to a non-zero eigenvalue $Hu=\lambda u$, then 
$v=Du$ is an eigenvector to the same eigenvalue $Hv = \lambda v$. But if
$u$ is an even dimensional form, then $v$ is an odd dimensional form. We
have now an isomorphism between even and odd parts of the image of $H$. 
It follows that the super trace of $H^k$ is zero for $k>0$. 
As ${\rm str}(H^0)={\rm str}(1)=\chi(G)$, when we apply the super trace to 
$e^{-t H} = \sum_{k=0}^{\infty} t^k H^k/k!$ this number is independent of $t$
and equal to the super trace of $H^0$ which is $\chi(G)$. 
\end{proof} 

\paragraph{}
The following result is the {\bf Poincar\'e-Hopf theorem}. 
\index{Poincar\'e-Hopf theorem}

\satz{ $\chi(G) = \sum_{k=0}^q (-1)^k f_k =  \sum_{k=0}^q (-1)^k b_k$.  }
\begin{proof}
Use the heat flow $e^{-t H}$ and the McKean-Singer symmetry. 
For $t=0$, the super trace of $e^{-t H}$ is the left hand side. 
For $t=\infty$, we get the right hand side. Since the super trace does not 
change, the left and right hand side agree. 
\end{proof} 

\paragraph{}
The {\bf Hodge decomposition} tells that any $k$-form $g$ can be written as a sum of an exact, a co-exact and a 
harmonic $k$-form: $g= df + d^* h + k$ and that harmonic forms are the intersection of the kernel of $d$
and $d^*$. 
\index{Hodge decomposition}

\satz{ ${\rm im}(H) = {\rm im}(d) + {\rm im}(d^*)$, ${\rm ker}(H)={\rm ker}(d) \cap {\rm ker}(d^*)$.  }

\begin{proof}
The operator $H: \Lambda \to \Lambda$ is symmetric so that image and kernel are
perpendicular. The image of $H$ splits into two orthogonal components ${\rm im}(d)$ and ${\rm im}(d^*)$.
That $f=dg$ and $h=d^*k$ are orthogonal is seen from 
$\langle dg,d^*k \rangle = \langle d dg,k \rangle = 0$. \\
To the second statement: the kernel of $H$ is contained both in the kernel of $d$ 
and the kernel of $d*$ because $Hf=0$ implies 
$0=\langle f,Hf \rangle=\langle d^*f,d* f \rangle + \langle df,df \rangle = |d^*f|^2+|df|^2$,
so that both $d^* f=0$ and $df=0$. On the other hand, if $df=0$ and $d^*f=0$, 
then $Hf=d d^*f + d^*d f=0$. 
\end{proof}

\paragraph{}
Let $d_k: \Lambda_k \to \Lambda_{k+1}$ be the restriction of $d$
to $\Lambda_k$. If the kernel ${\rm ker}(d_k)$ of dimension $z_k$ and the range 
${\rm im}(d)$ has dimension $r_k$, the rank-nullity theorem in linear algebra tells
$\dim(\ker(d_k)) + \dim({\rm im}(d_k))=f_k$, we get $z_k = f_k-r_k$.

\satz{$b_k={\rm dim}({\rm ker}(d_k))/{\rm im}(d_{k-1})$ and so $b_k=z_k-r_{k-1}$.}
\begin{proof}
${\rm im}(d_{k-1})$ is contained in ${\rm ker}(d_k)$. The $n \times f_k$ matrix
$D_k$ containing the blocks $d_{k-1}^*$ and $d_k$ satisfies $D_k^* D_k = H_k$. 
Now use that ${\rm im}(d_{k-1})$ is the orthogonal complement of ${\rm ker}(d_{k-1}^*)$
and row reduce $F_k$ to see that $b_k$, the number of redundant columns in $F_k$
are the columns, where no leading $1$ for $d_{k-1}^*$ and no leading $1$ for $d_k$ appears.
\end{proof}

\paragraph{}
Adding $z_k = f_k-r_k$ and $b_k=z_k-r_{k-1}$ gives  $f_k-b_k = r_{k-1}+r_k$. 
Taking the alternating sum on both sides using $r_{-1}=0$ and $r_k=0$ for $k>q$, 
telescopes $\sum_{k=0}^{q} (-1)^k (f_k-b_k) = \sum_{k=0}^{\infty} (-1)^k (r_{k-1}+r_k) = 0$.
This is the traditional proof of the Euler-Poincar\'e theorem. 

\paragraph{}

{\bf Examples}.
1) $b_0$ is the number of components, $b_1$ the number 1 dim holes. 
2) For connected orientable surfaces, we have $b_0=b_2=1$ so that 
      $\chi(G)=|V|-|E|+|F| = 2-b_1 =2-2g$ where $g$ is the genus of the surface. 
3) For connected regions in the plane we have $|V|-|E|+|F|=1-g = b_0-b_1$.
4) A contractible graph has $b_0=1$ and $b_k=0$ for all $k>0$. 
5) For a q-sphere $b_0=b_q=1$ and $b_k=0$ for all $0<k<q$. 
 \vfill \pagebreak
\begin{center}\scalebox{0.8}{\includegraphics{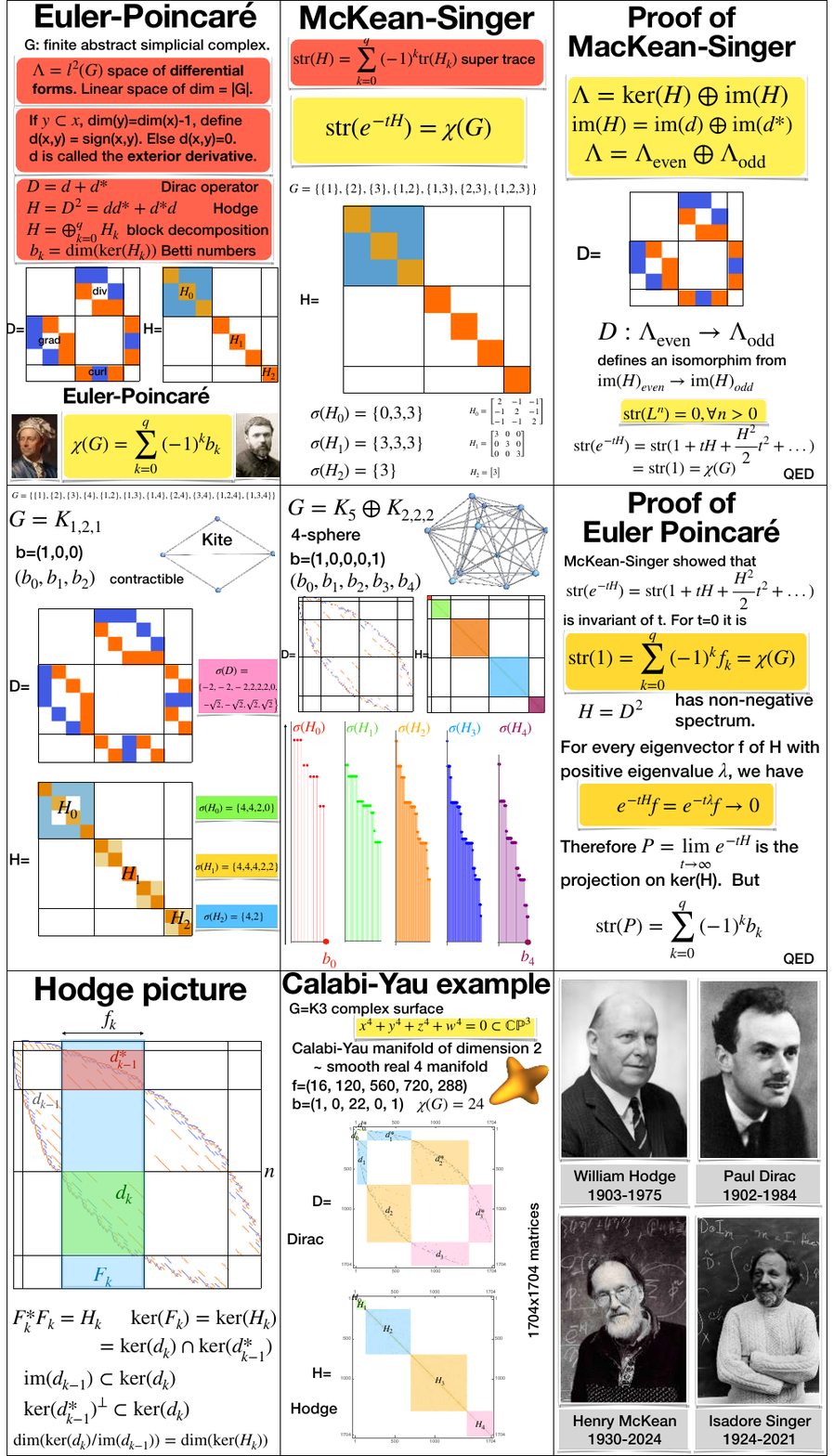}}\end{center} \pagebreak
\setcounter{section}{6}  \chapter{Unit 6: Unimodularity}

\summary{
A simplicial complex $G$ defines a connection matrix $L$ that
is unimodular. The inverse matrix $g$ of $L$ is explicitly given
in the form of the Green-Star formula 
$g(x,y)=\omega(x) \omega(y) \chi(U(x) \cap U(y))$. 
The total potential $\sum_{x,y} g(x,y)$ is $\chi(G)$. 
}

\paragraph{}
If $G$ is a simplicial complex,
the {\bf stars} $U(x)=\{ y \in G, x \subset y\}$ are open sets,
while the {\bf cores} $C(x) \{ y \in G, y \subset x \}$ are closed. 
For $A \subset G$, the {\bf Fermi characteristic}
$\phi(G)=\prod_{x \in G} \omega(x)$ of is a multiplicative sibling of
the additive Euler characteristic 
$\chi(A)=\sum_{x \in A} \omega(x)$, where $\omega(x)=(-1)^{{\rm dim}(x)}$.
Define the {\bf connection matrix} $L(x,y) = \chi(C(x) \cap C(y))$. The entries
are  $1$ if $x \cap y \neq \emptyset$ and $0$ else. 
Define also $\psi(G) = {\rm det}(L(G))$. The {\bf unimodularity theorem} is:
\index{star}
\index{core}
\index{Fermi characteristic}
\index{open set}
\index{closed set}
\index{unimodularity theorem}

\satz{ ${\rm det}(L) = \phi(G) \in \{-1,1\}$.  }

\begin{proof}
The simplicial complex $G$ can be built up recursively from a smaller complex
$H$ by adjoining a new $k$-simplex $x$ to an already existing $(k-1)$-sphere $A$
in $H$ that is not yet the boundary of a ball in $H$. This enlargement $H \to G=H \cup_A \{x\}$
changes the Fermi characteristic $\phi(G)$ to $\phi(G) = \phi(G) \omega(x)$.
A path interpretation of the determinant verifies that this also
changes the determinant $\psi(G) = {\rm det}(L)$ by a factor $\omega(x)$.
It is the content of the  multiplicative Poincar\'e-Hopf analog below.
\end{proof}

\paragraph{}
{\bf Poincar\'e-Hopf} is the additive statement $\chi(G \cup_A \{x\})= \chi(G) + (1-\chi(A))$.
It used the inclusion-exclusion formula and $\chi(B(\{x\}))=1$.  Multiplicatively, we have:
\index{Multiplicative Poincare-Hopf}
\index{inclusion-exclusion}

\satz{ $\psi(G \cup_A \{x\}) = \psi(G) (1-\chi(A))$.}
\begin{proof}
$Y: A \to (\psi(G \cup_A x) - \psi(G))$ is a valuation that
satisfies $\psi(G \cup_A x)=0$ if $A$ is a complete complex so that
$Y(A)=-\psi(G)$ if $A$ is a complete complex.
The scaled valuation $-Y(A)/\psi(G)$ takes the value $1$ for complete sub-complexes.
This characterizes Euler characteristic: we have
$-Y(A)/\psi(G)=\chi(A)$ meaning $Y(A) = -\chi(A) \psi(G)$
But $(\psi(G \cup_A x) - \psi(G)) = -\chi(A) \psi(G)$ is equivalent to 
$\psi(G \cup_A \{x\}) = \psi(G) (1-\chi(A))$.
\end{proof} 

\paragraph{}
Let $g(x,y) = \omega(x) \omega(y) \chi(U(x) \cap U(y))$ be the 
{\bf Green function}. It is like $L$ a $n \times n$ matrix. 
It invokes the {\bf stars} $U(x),U(y)$ of $x$ and $y$. 
The next result is the {\bf Green-star formula}:
\index{Green-Star Formula}

\satz{ $L^{-1} = g$.  }

\begin{proof}
Write the matrix entries of $L$ and $g$ are linear expressions in $\omega(x)$:
$$  L(u,v) = \sum_{x \in G, x \subset u \cap v} \omega(x)                      \; , $$
$$  g(v,w) = \sum_{y \in G, v \cup w \subset y} \omega(v) \omega(w) \omega(y)  \; . $$
Now multiply the two matrices $\sum_v L(u,v) g(v,w)$. This is a sum over $x,y,v$. 
Distinguish two cases. If $u=w$, we divide the sum up into two parts, one part where
$x =y$ which implies $u=w,x,y=v$ and giving $1$ and a part 
where $x \neq y$, where the sum over $v$ is zero. The diagonal entries are $1$. \\
In the second case $u \neq w$ the situation is similar, but $x=y$ is no more possible. 
We end up with 0.
\end{proof}

\paragraph{}
If $g(x,y)$ is interpreted as the {\bf potential energy} between the
simplices $x$ and $y$, the number $E(G)=\sum_{x,y \in G} g(x,y)$ is the
{\bf total energy}. The {\bf energy theorem} is:
\index{potential energy}
\index{total energy}

\satz{ $\sum_{x,y \in G} g(x,y) = \chi(G)$ }

\begin{proof} 
We have $g(x,x)=\chi(U(x)) = 1-\chi(S(x))$. Call it $i(x)$.
By the Cramer formula, $g(x,x)$ is $\psi(G \setminus x)/\psi(G)$ 
by the multiplicative Poincar\'e-Hopf theorem this is $1-\chi(S(x))$.
Now think of $G$ is the vertex set of a graph in which two vertices are
connected if one is contained in the other. Take the functional ${\rm dim}(x)$. 
Hyperbolicity gives $S(x)=S^-_f(x) \oplus S^+_f(x)$ so that
$i_{f}(x) i_{-f}(x) = i(x)$. Since $i_f(x) = \omega(x)$ and $i_{-f}(x)= 1-\chi(U(x))$
$k(x) = \omega(x) (1-\chi(S(x)))$ is also a curvature adding up to $\chi(G)$. 
This implies $\sum_x \omega(x) \chi(S(x))=0$ which can be read as a
McKean-Singer  ${\rm str}(L^k)=\chi(G)$ for $k=-1,0,1$.
Introduce the potential $V(x) = \sum_{y} g(x,y)$. It satisfies
$V(x)= \omega(x) g(x,x) = k(x)$. Gauss-Bonnet $\sum_x k(x)=\chi(G)$ now proves
the theorem. 
\end{proof} 
\index{energy theorem}

\paragraph{}
Examples: \\
1) ${\rm det}(L(K_1))=1$ but ${\rm det}(L(K_k))=-1$ for $k>1$. \\
2) The diamond complex
$G = \{$ $(1),(2)$, $(3),(4)$, $(1,2),(1,3)$, $(2,3),(2,4)$, $(3,4)$, $(1,2,3),(2,3,4)\}$
has  $f$-vector $f=-(4,5,2)$, Euler characteristic $\chi(G)=4-5+2=1$ and Betti vector $b=(1,0,0)$.
Its Fermi characteristic is $(+1)^4 (-1)^5 (+1)^2=-1$. 
\index{diamond complex}

 \vfill \pagebreak
\begin{center}\scalebox{0.8}{\includegraphics{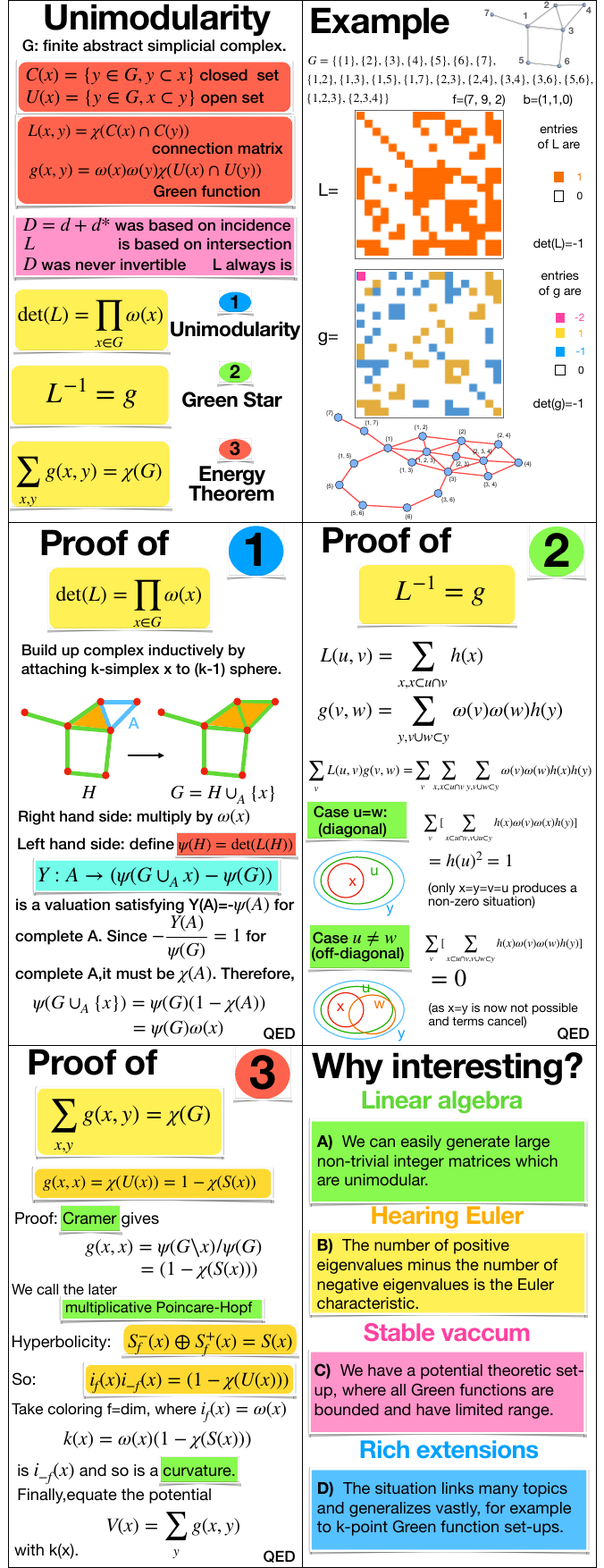}}\end{center} \pagebreak
\setcounter{section}{7}  \chapter{Unit 7: Brouwer Lefschetz}

\summary{
The sum of the indices of a simplicial map $T$ on a simplicial complex
$G$ is equal to the Lefschetz number $\chi_T(G)$, the super trace of
the Koopman operator of $T$ restricted to the harmonic forms of $G$.  }

\paragraph{}
A map $T: G \to G$ on a finite abstract simplicial complex $G$ 
is a {\bf simplicial map} if it maps the set $V$ of $0$-dimensional simplices 
into itself and if it is order-preserving, meaning that $x \subset y$ implies $T(x) \subset T(y)$. 
The map $T$ does not have to be invertible. A constant map $T$ that maps $G$ to $x \in G$ 
for example is a simplicial map only if $x$ is $0$-dimensional.
\index{simplicial map}
\index{order preserving map}

\paragraph{}
If $T$ is a simplicial map, its {\bf attractor} $\bigcap_{k \geq 0} T^k(G)$ is 
a simplicial complex on which $T$ is invertible. As we are interested in fixed 
points of $T$ and fixed points must be in the attractor, we can assume $T$ is invertible.
\index{attractor}

\paragraph{}
An invertible simplicial map $T$ induces a {\bf Koopman operator} $U$ on the Hilbert space 
$\Lambda$ of differential forms. It is $U f(x) = f(Tx) {\rm sign}(T|x)$. 
If $T$ is invertible, it is a unitary operator. 
The linear map $U$ commutes with the exterior derivative $d$ and so 
commutes with the Hodge Laplacian. It therefore preserves the linear space of 
$k$-forms $\Lambda_k$ and also preserves the harmonic $k$-forms. 
\index{Hodge Laplacian}
\index{Koopman operator}

\paragraph{}
Let $U_k$ be the restriction of $U$ to $k$-forms. Its super trace
${\rm str}(U|{\rm ker}(H)) = \sum_{k=0}^q (-1)^k {\rm tr}(U_k|{\rm ker}(H_k))$ 
of $U$ is the {\bf Lefschetz number} of $T$. We call it $\chi_T(G)$ because if
$T$ is the identity, then $\chi_T(G)$ is the Euler characteristic $\chi(G)$. 
\index{super trace}
\index{Lefschetz number}

\paragraph{}
Let $F=\{ x \in G, T(x)=x \}$ denote the {\bf fixed point set}.
$T$ induces a permutation of the vertices of each fixed point $x$ for which ${\rm sign}(T|x)$ 
is the sign of. The number $i_T(x) = \omega(x) {\rm sign}(T|x)$ is 
called the {\bf index} of $x$. 
\index{fixed point set}
\index{index of a fixed point}

\paragraph{}
The Lefschetz fixed point theorem is

\satz{ $\chi_T(G) = \sum_{x \in F} i_T(x)$.  }
\index{Lefschetz fixed point theorem}

\begin{proof}
The operator $U$ commutes with the heat flow $e^{-t L}$. Define the 
operator $U_t = U e^{-t L}$. The McKean-Singer symmetry implies that
${\rm str}(U_t)$ is independent of $t$. For $t=0$, it is the super trace of 
$U$ which is $\sum_{x \in F} i_T(x)$ by the theorem below.
Since $\lim_{t \to \infty} e^{-t L} = P$ is the projection on the harmonic
forms, we have ${\rm str}(U e^{-t L}) \to \chi_T(G)$, proving the theorem. 
\end{proof} 

\satz{ ${\rm str}(U) = \sum_{x \in F} i_T(x)$ }

\begin{proof}
$\sum_{x \in G} \omega(x) U(x,x) = \sum_{x \in F} \omega(x) U(x,x) 
  = \sum_{x \in F} \omega(x) {\rm sign}(T|x) = \sum_{x \in F} i_T(x)$. 
\end{proof} 

\paragraph{}
Lets see what happens if $G$ is the complex of a complete graph $K_n$
and $T$ is a permutation of $V=\{1, \dots, n\}$. In this case ${\rm str}(U)=1$
which is $\chi_T(G)$. The map $U$ encodes the fixed points in the diagonal
If the permutation belonging to $T$ is $\pi = (c_1)(c_2) \dots (c_m)$, where $c_i$ are the
cycles, then every non-empty collection of cycles produces a fixed simplex of $T$.

\paragraph{}
The Lefschetz theorem for contractible spaces is {\bf Brouwer's fixed point theorem}:
\index{Brouwer's fixed point theorem}

\satz{A map on a contractible $G$ has at least $1$ fixed point.}
\begin{proof}
If $G$ is contractible, the cohomology only consists of
constant $0$-forms. The Lefschetz number is $1$.
\end{proof}

\paragraph{}
Let ${\rm Aut}(G)$ denote the {\bf automorphism group} of $G$, the set of
symmetries. Given a subgroup $A$ of ${\rm Aut}(G)$, the {\bf average Lefschetz number} is
$\chi(G/A) = \frac{1}{|A|} \sum_{T \in A} \chi_T(G)$. If $G$ is trivial it is $\chi(G)$.
\index{automorphism group}
\index{average Lefschetz number}
\index{quotient complex}

\satz{If $H=G/A$ is a complex, $\frac{1}{|A|} \sum_{T \in A} \chi_T(G) = \chi(G/A)$.}

\begin{proof}
Burnside's theorem for $A$ acting on $G$ gives $f_k(G/A) = \frac{1}{|A|} \sum_{T \in A} |F_k(T)|$
with $F_k(T)=\{ x \in G_k, T(x)=x\}$ so that
$\chi(G/A) = \sum_{k=0}^{q} (-1)^k f_k(G/A)$ is
             $$\frac{1}{|A|} \sum_{T \in A} (-1)^k |F_k(T)|
            = \frac{1}{|A|} \sum_{T \in A} \sum_{x \in F} i_T(x) 
            = \frac{1}{|A|} \sum_{T \in A} \chi_T(G) \; . $$
\end{proof}

\paragraph{}
The {\bf curvatures} of $G/A$ are $\kappa(x) = \frac{1}{|A|} \sum_{T \in A} i_T(x)$. They add up to $\chi(G/A)$. 
We can see $G$ as {\bf branched cover} over $H=G/A$. If $A$ is the cyclic group generated by $T$, then $G$ is a 
$k:1$ cover over $H=G/A$. The fixed points of $T$ are branch points. If $T,T^2,\dots, T^{k-1}$ have 
no fixed points then $\chi_{T^l}(G)=0$ and $\chi_{Id}(G)=\chi(G)$ so that 
$\chi(H)=\chi(G)/|A|$, a special {\bf Riemann-Hurwitz} result.
\index{Riemann-Hurwitz}
\index{branched cover}
\index{curvatures of a quotient}

\paragraph{}
{\bf Examples:} 
1) For $G=[0,1,2,3]$ and $T(x)=T(3-x)$, there is single fixed point $(1,2)$ of index $1$.
2) If $G$ a figure $8$ complex that carries an involution with 3 fixed points of index $1$. 
The trace on harmonic $0$ forms is $1$. The trace on harmonic $1$ forms is $-2$.  
3) Invertible orientation preserving maps on a $2k$-sphere,
or invertible orientation reversing  maps on a $(2k+1)$-sphere have at least $2$ fixed points. 

 \vfill \pagebreak
\begin{center}\scalebox{0.8}{\includegraphics{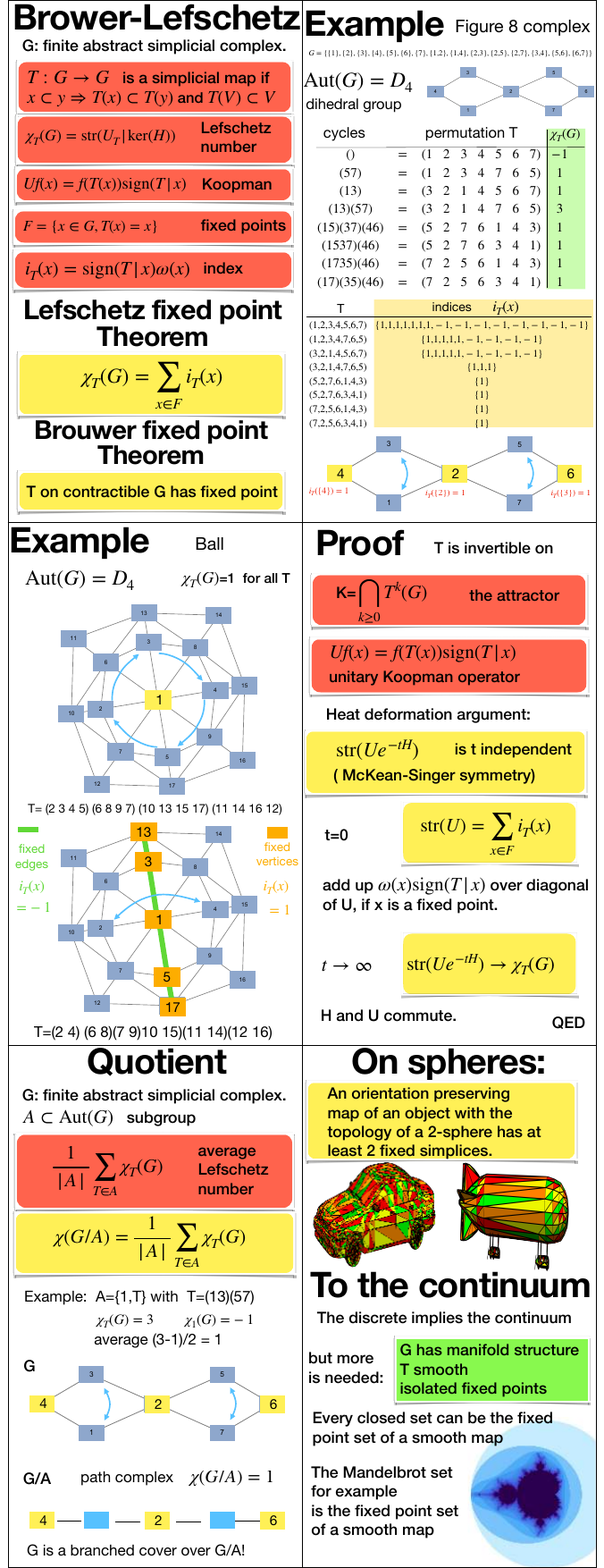}}\end{center} \pagebreak
\setcounter{section}{8}  \chapter{Unit 8: Sphere formula }

\summary{
If $G$ is any finite abstract simplicial complex, 
then $\sum_{x \in G} \omega(x) \chi(S(x))=0$. 
If all unit spheres have the same non-zero 
Euler characteristic, like for odd-dimensional manifolds, 
then the Euler characteristic of the complex is zero.  }

\paragraph{}
Let $G$ be a finite abstract simplicial complex. 
The {\bf unit sphere} $S(x)$ of $x \in G$ is the boundary of 
the {\bf star} $U(x)  = \{ y \in G, x \subset y\}$, the smallest
open set that contains $x$. As this 
is $\overline{U}(x) \setminus U(x)$, the unit sphere
is a simplicial complex.  The topological definition links
it to a natural topology $\mathcal{O}$ on $G$, even so it 
is never Hausdorff, if the dimension of $G$ is positive. 
\index{unit sphere in simplicial complexes}
\index{star}
\index{open set}
\index{non-Hausdorff}

\paragraph{}
For a vertex $v$ in a graph, the unit sphere $S(v)$ is the graph generated by 
the neighbors of $v$. For a simplex $x$ in a complex $G$, the unit sphere is $S(x)=\delta U(x)$. 
There is a connection: if $\Gamma(G)=(V,E)$ is the graph
with vertices $V=G$ and $E=\{(x,y), x \subset y \; {\rm or} \; y \subset x \}$,
one can associate $S(x)$ with $\Gamma(S^+(x)) \oplus \Gamma(S^-(x))$. 
We write $S(x) = S^+(x) \oplus S^-(x)$ even-so we join an open with a closed set. 
This is a form of hyperbolicity because for a hyperbolic fixed point of a 
smooth map $T$ in $\mathbb{R}^n$, we have
{\bf stable and unstable manifolds} $W^{-}(x)$ and $W^+(x)$ which when intersected with a small
sphere $S(x)=S_r(x)$ produce two spheres $S^-(x)$ and $S^+(x)$ which have the property that
their join is $S(x)$. In $G$, the hyperbolic structure comes from the
{\bf dimension functional} ${\rm dim}$ on $G$. Here is the {\bf hyperbolicity lemma}:
\index{hyperbolicity lemma}

\satz{$S(x) = S^+(x) \oplus S^-(x)$ for all $x \in G$. }

\begin{proof} 
We show the graph join $\Gamma(S(x)) = \Gamma(S^+(x)) \oplus \Gamma(S^-(x))$. A simplex
$y \in G$ is a vertex in $\Gamma(S(x))$ if either $y$ is a subset of some $z \in S(x))$
meaning $y$ is in $S^-(x)$ or if $y$ contains some $z \in S(x)$ meaning
$y$ is in $S^+(x)$. So, the vertex sets of the two graphs agree. If two points
$y,z$ in $\Gamma(S(x))$ are connected, we either have $y \subset z$ or $z \subset y$. 
There are four cases: either both are in $S^-(x)$ which means they are connected 
in $S^-(x)$ and so connected in the join. The same hold if both are in $S^+(x)$. 
Any pair $y \in S^-(x)$ and $z \in S^+(x)$ satisfies $y \subset x \subset z$ 
so that they are connected in $S(x)$. 
\end{proof}


\satz{ $\sum_{x \in G} \omega(x) \chi(S(x))=0$.  }
\index{sphere formula}

\begin{proof}
(i) For any $x$ we have the local valuation formula 
$\chi(B(x))=\chi(U(x))+\chi(S(x))$.
Proof: For any subsets $A,B \subset G$, whether open or closed or neither, we have
the valuation formula $\chi(A) + \chi(B) = \chi(A \cup B) - \chi(A \cap B)$ because
each of the basic valuations $f_k(G)$ counting the number of $k$-dimensional
simplices satisfies the formula and $\chi$ is a linear combination of such
basic valuations.
(ii) $\chi(B(x))=1$ for all $x$.
Proof: Use induction with respect to the number of elements in $B(x)$.
If $B(x)$ has one element, it has $\chi(B(x))=1$. We can reduce the size of $B(x)$ by
taking a way an element $y \in B(x)$ different from $x$. 
The complex $B'(x)=B(x) \setminus U(y)$ is
now a unit ball $B'(x)$ in a smaller $G \setminus U(y)$.  The induction assumption assures
that $\chi(B'(x))=1$. The local valuation formula gives that $\chi(B(x)))=\chi(B'(x)$. \\
(iii) the energy formula $\sum_{x \in G} \omega(x) \chi(U(x)) = \chi(G)$ was proven 
in the unimodularity unit. It told that the super trace of $g$ was Euler characteristic. 
Here is an other proof: Replace the $\pm 1$-valued function 
$\omega(x)$ by a more general function $h(x)$ and define
$\chi_h(A) = \sum_{x \in A} h(x)$ and $\psi_h(x) = \omega(x) \chi_h(U(x))$.
Both maps $h \to \chi_h(A)$ and $h \to \psi_h(A)$ are linear in $h$ for fixed $A$
so that we only need to consider the case where $h(x_0)=1$ and $h(x)=0$ for $x \neq x_0$.  
The right hand side is then $\chi_h(G) =1$. 
The left hand side is $\sum_{x \in G} \omega(x) \chi_h(U(x)) =
\sum_{x \in G, x \subset x_0} \omega(x)$, which is the Euler characteristic 
of the simplicial complex $\overline{ \{x_0\} } = \{ x \subset x_0 \}$ which is always $1$. 
(iv) $\omega(G) = \sum_{x \in G} \omega(x) \chi(B(x)) = 
   \sum_{x \in G} \omega(x) = \chi(G)$. 
\end{proof} 

\paragraph{}
The value $g(x,x)=\omega(x)^2 \chi(U(x) \cap U(x)= \chi(U(x))$ 
is a diagonal entry in the {\bf Green function}. It is the 
{\bf self-energy} of $x \in G$. Now $\chi(S(x)) = 1-\chi(U(x)$ gives
$\sum_{x \in G} \omega(x) \chi(S(x)) = \sum_{x \in G} \omega(x) (1-\chi(U(x)) = \chi(G) - \chi(G) = 0$, 
which is the sphere formula.  

\paragraph{}
What is $\sum_{x \in G} \chi(S(x))$? This is the trace of $L-g$. We called this the 
{\bf Hydrogen functional} because of formal similarity with the quantum mechanical
Hamiltonian of the Hydrogen atom, where the potential part involves the kernel of
the inverse of the kinetic part - the $1/r$ potential being determined by 
the Laplacian $\Delta$ in $\mathbb{R}^3$ for example. We called $L-L^{-1} = L-g$ the {\bf Hydrogen
operator}. 
\index{hydrogen operator}
\index{hydrogen formula}

\satz{ $\sum_{x \in G} \chi(S(x)) = {\rm tr}(L-L^{-1})$.}

\paragraph{}
The Euler-Gem and sphere formula together immediately give the zero
Euler characteristic result:

\satz{ All odd-dimensional manifolds satisfy $\chi(G) =0$.  }

\begin{proof}
Euler-Gem formula assures that odd-dimensional
manifolds have unit spheres with $\chi(S(x))=2$. The sphere formula
implies $0=\sum_x \omega(x) \chi(S(x)) = \sum_x \omega(x) 2 = 2 \chi(G)$.
\end{proof}

 \vfill \pagebreak
\begin{center}\scalebox{0.8}{\includegraphics{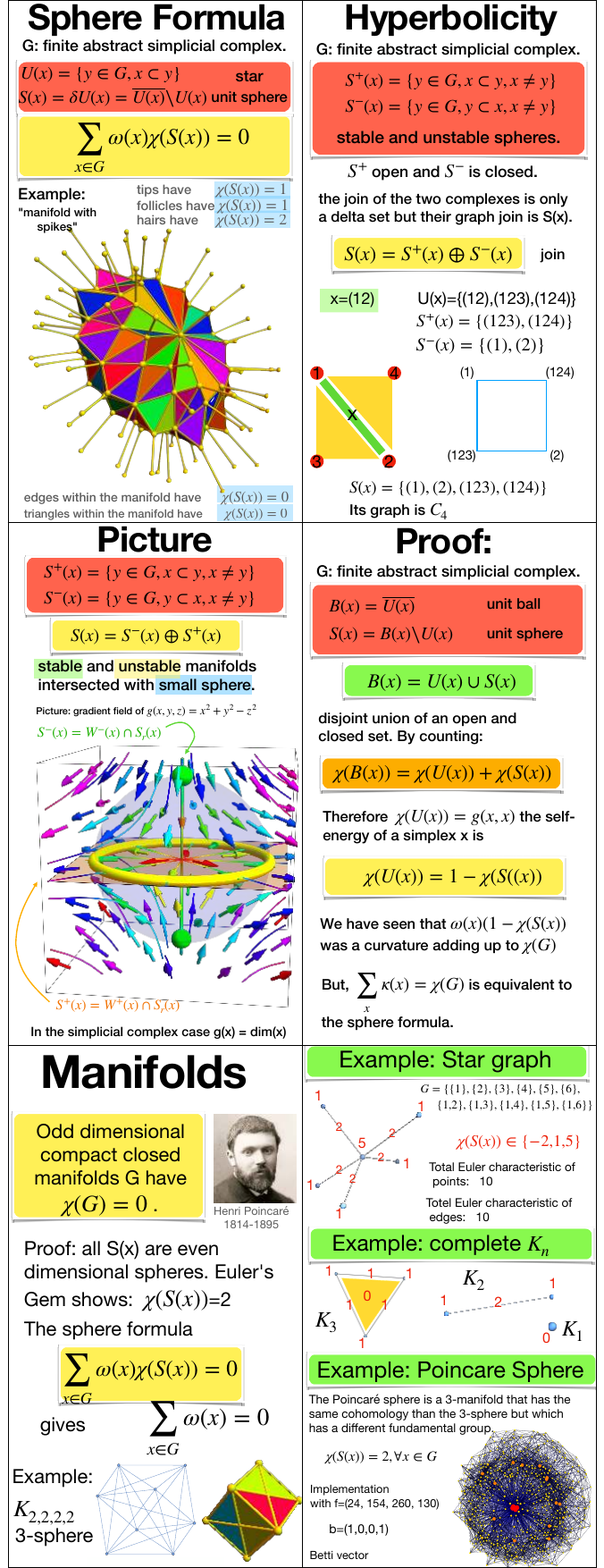}}\end{center} \pagebreak
\setcounter{section}{9}  \chapter{Unit 9: Level sets }

\summary{
If $G$ is a discrete q-manifold and $g: V \to K_k=\{0, \dots, k\}$ is a function,
then $M_g = \{ x \in G | g(x)=K_k \}$ is a $(q-k)$-manifold if not empty.  }

\paragraph{}
A {\bf $q$-manifold} is a finite abstract simplicial complex $G$ such that 
every unit sphere $S(x)$ is a $(q-1)$-sphere. A {\bf $q$-sphere} is a
$q$-manifold $G$ such that $G \setminus U(x)$ and $S(x)$ are both 
contractible. A complex $G$ is {\bf contractible} if there exists
$x$ such that both $G \setminus U(x)$ and $S(x)$ are both contractible. 
The initial object {\bf void} $0$ is the $(-1)$ sphere and the terminal
object $1$ is contractible. 
\index{q-manifold}
\index{contractible complex}

\paragraph{}
The vertex set $V=\bigcup_{x \in G} x$ in $G$ can be identified with 
the $0$-dimensional objects $G_0$ in $G$. A function $g: V \to K_k=\{0, \dots, k\}$ 
defines a {\bf level set} of $G$ denoted by $M_g = \{ x \in G, g(x) = K_k\}$, where 
$g((x_0,x_1,\dots, x_))=\{ g(x_0), \dots g(x_k)\}$
In the case $k=1$ for example, this is the set of $x$ on which $g(x)=\{0,1\}$
takes both values, meaning that $\{ v \in x, g(v)=0\}$  and $\{v \in x, g(v)=1\}$ 
are both non-empty. The set $M_g$ is a priori only open. This means 
$x \in M_g$, and $y \in G$ contains $x$ then $y \in M_g$. But $M_g$ defines a 
simplicial complex: take the set $M_g$ as the vertices of a graph $\Gamma$ and connect 
two vertices $(a,b)$ if $a \subset b$ or $b \subset a$. The Whitney complex of this graph is 
then a simplicial complex and if we say $M_g$ is a manifold, we mean that its simplicial
complex is. 
\index{level set}

\paragraph{}
\satz{Any $g:V(G) \to K_k$ on a q-manifold $G$ has a level set $M_g$ that is a
$q-k$ manifold if it is not empty. }
\index{level set theorem}
\index{submanifold theorem}

\paragraph{}
The proof uses the hyperbolic structure of 
any simplicial complex. In the case of manifolds, we can even deal with a
Morse complex because every $x \in G$ is a critical point of $x \to {\rm dim}(x)$. The 
Morse index of $x$ is $m(x)={\rm dim}(S^-(x)) = k$, where $S^-(x)=\{y \in G, y \subset x, y \neq x \}$.
Complementary is the $S^+(x) = \{ y \in G, x \subset y, y \neq x \}$. 
$G$ is a manifold if and only if every $S^+(x)$ is a sphere. Again, $S+(x)$ is only an open set
at first and we refer to the Whitney complex of its graph. Hyperbolicity is:
\index{hyperbolicity}

\satz{If $G$ is a $q$-manifold and ${\rm dim}(x)=k$, then $S(x)$ is the join of 
the $(k-1)$-sphere $S^-(x)$ and $(q-k-1)$ sphere $S^+(x)$.}

\begin{proof}
The graph of $S(x)$ contains both simplices $S^-(x)$ strictly contained in $x$ and simplices
$S^+(x)$ that strictly contain $x$. Every $y \in S^-(x)$ is contained in any $z \in S^+(x)$. 
\end{proof} 

\paragraph{}
In the special case, when $G$ is the boundary complex $S^-(x)$ of 
a $k$-simplex $x$, level sets are always spheres. 
For the 2-simplex $x=K_3$ for example, the boundary sphere is
$S^-(x)=\{ \{1\},\{2\},\{3\},\{1,2\},\{2,3\},\{1,3\}\}$.
For any function $g: \{1,2,3\} \to \{0,1\}$ the level set is always a $0$-sphere. 

\paragraph{}
A special case is if $G$ be a complete $n$-dimensional complex. Every surjective
$g:V(G) \to K_k$ produces a level set $M_g$ that is a $(n-k)$-sphere. 
The sets $V_j = \{ g=j \} \subset V$ partition $V$ into $k+1$ sets of size $n_j>0$ such that the sum is $n+1$.
The level set $M_g$ is the join of all these $(n_j-1)$-spheres which is a $(n-k)$ - sphere
because $(\sum_j n_j)-1= n+1- k-1 = n-k$. 

\paragraph{}
Now the proof of the main theorem:

\begin{proof}
We use induction with respect to dimension $q$.
The unit sphere $S_G(x)$ in $G$ is the join of $S^+(x)$ and $S^-(x)$. 
The unit sphere $S_M(x)$ in $G$ is the join of $S_M^+(x)$ and $S_M^-(x)$. 
We have $S_M^+(x)=S^+(x)$ because if $y$ is in $S(x)$ and $g$ takes all values
on $x$, then $g$ takes all values on $y$ so that $y$ is also in $S_M^+(x)$. 
As seen above, the level surface in the $k-1$ sphere $S^-(x)$ is $k-2$-sphere. 
So $S_M(x)$ is the join of a $(k-2)$ sphere with a $q-k-1$ sphere which is a
$q-k-1$ sphere
\end{proof}

\paragraph{}
The set-up can be generalized if $K_k$ is replaced with a $(k+1)$-partition complex
$P=K_{n_0,n_1, \dots n_k}$. If a function $g:V(G) \to V(P)$ is given, define
$M_g = \{ x \in G, g(x)$ contains at least one facet of $P \}$. In that case again
$M_g$ is a $q-k$ manifold if it is not empty. Every facet $f$ of $P$ now produces
a patch $M_{g,f} = \{ x \in g, f \subset g(x) \}$ which is a $q-k$ manifold with 
boundary. The union of all these patches is the manifold $M_g$ without boundary. 

\satz{Any $g:V(G) \to P$ from a q-manifold $G$ to a {\bf partition complex} $P$ of dimension $k$
has a level set $M_g$ that is a $q-k$ manifold or $\emptyset$. }

\paragraph{}
In the continuum, level sets in $q$-manifolds are not manifolds in general, like
for $g(x,y) = \sin(x) \sin(y) = 0$ on the 2-manifold $M =\mathbb{T}^2$. The Morse Sard theorem 
assures that for almost all $c$ in the range of $g: M \to \mathbb{R}^k$, 
the level set $\{ g=c\}$ is a $q-k$ manifold. In the discrete, no regularity at all is needed. 
\index{Morse-Sard theorem}

\paragraph{}
{\bf Examples:}
\begin{itemize}
\item If $g$ is $1$ on some vertex $v$ and $0$ else, then 
$M_g = S(v)$ is the unit sphere.
\item If $G$ is a hexagonal lattice and $g: V \to \{0,1\}$ is a function, 
then $M_g$ consists of all triangles and edges on which $g$ takes two values. 
If an edge $e$ is in $M_g$ then the two attached triangles are there. 
If $t$ is a triangle in $M_g$, then there are exactly two edges contained
in $M_g$ where $g$ takes two values. Since every point in $M_g$ has exactly 
two neighbors, the level set must be a union of circular graphs. 
\end{itemize}

 \vfill \pagebreak
\begin{center}\scalebox{0.8}{\includegraphics{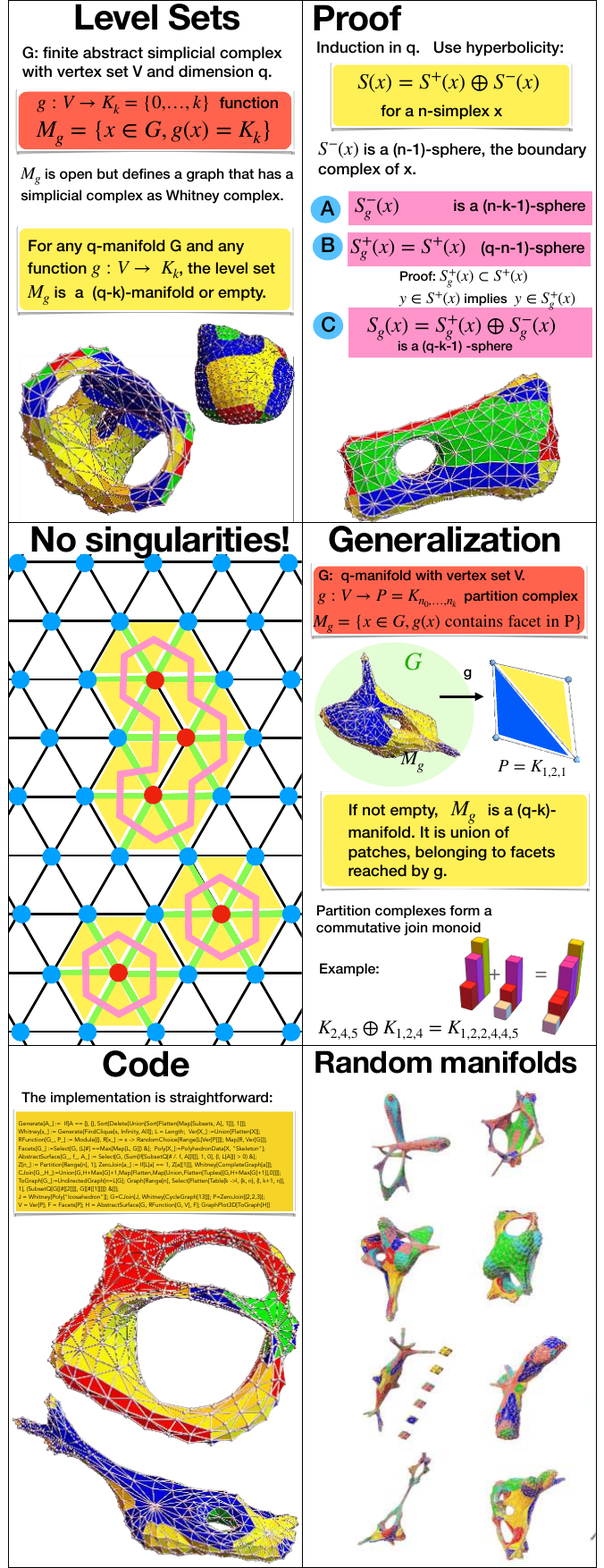}}\end{center} \pagebreak
\setcounter{section}{10} \chapter{Unit 10: Index formula}

\summary{
If $g$ is a vertex coloring of a graph $(V,E)$, the level set $S_g(v) \subset S(v)$
produces the {\bf symmetric Poincare-Hopf index} $j_g(v) = 1-\chi(S(v))/2 - \chi(S_g(v))/2$.
This implies that odd-dimensional manifolds are flat and that for 4-manifolds,
curvature $K(v)$ is the genus expectation of a random 2-manifold in the 3-manifold $S(v)$. }

\paragraph{}
If a finite abstract simplicial complex $G$ defines a graph $(V,E)$ the Poincar\'e-Hopf theorem
for $G$ is just the Poincar\'e-Hopf theorem for the corresponding graph with $V=G$. 
A locally injective function $g:G \to R$ is a vertex coloring function $g: V \to R$. The index 
$i_g(v)=1-\chi(S^-_g(v))$ defines then the {\bf symmetric index} $j_g(v)=[i_g(v)+i_{-g}(v)]/2$.
Since ${\rm sign}(g(w)-g(v))$ is a $\{-1,1\}$-valued function on $S(v)$, we can look at its level set
$S_g(v)$ in $S(v)$. If $S(v)$ was a $(q-1)$-sphere, then $S_g(v)$ is a $(q-2)$-manifold. 
\index{symmetric index}
\index{manifold for the symmetric index}

\paragraph{}
The following result combines Poincar\'e-Hopf with the discrete Sard story. It shows that for
manifolds, the symmetric index is expressible as the Euler characteristic of sub-manifolds. 
First for general complexes:
\index{index formula}

\satz{If $g$ is a coloring, then $j_g(v) = 1-\chi(S(v))/2 - \chi(S_g(v))/2$. }

\begin{proof}
The level set result clarifies what we mean with $S_g(v)$: it is a level set in 
$S(v)$. We have a decomposition $S(v)=S_g^-(v) \cup S_g^+(v) \cup S_g(v)$ 
in the sense that every simplex in $S(v)$ is either a simplex in $S_g^-(v)$
or a simplex in $S_g^+(v)$ or a simplex in $S_g(v)$. Therefore,
$\chi(S(v)) = \chi(S_g^-(v)) + \chi(S_g^+(v)) - \chi(S_g(v))$, the negative sign coming
from the fact that $S_g(v)$ has edges as smallest dimensional simplices, shifting
dimensions by $1$. Therefore, $\chi(S(v)) + \chi(S_g(v)) = 1-i_g(v) + 1-i_{-g}(v) = 2-2 j_g(v)$. 
Now solve for $j_g(v)$. 
\end{proof} 

\paragraph{}
In the case of manifolds, we have by the Euler gem formula for odd $q$ that 
$j_g(v) = -\chi(S_g(v))/2$ which means that $j_g(v)=0$. Also the sphere formula
had shown that the odd-dimensional manifold $S_g(v)$ had zero
Euler characteristic. Therefore, for any index expectation curvature that comes
from a probability measure, for which $g \to -g$ is an isometry we have:

\satz{Odd-dimensional manifolds are flat. }

\paragraph{}
In the case of even-dimensional manifolds we have $j_g(v) = 1-\chi(S_g(v)))/2$. 
Motivated from notions in two dimensions, let us call $g=1-\chi(G)$ the {\bf genus} 
of $G$. It is also known as a {\bf reduced Euler characteristic}.
\index{genus of a manifold}
\index{reduced Euler characteristic}

\satz{Curvature $K(v)$ in even dimensions is the genus average of $q-2$ manifolds in $S(v)$.}

\paragraph{}
For a 4-manifold for example, where $S_g(v)$ is an orientable manifold of genus $g$
with Euler characteristic $\chi=2-2g$, we have:

\satz{For a 4-manifold, the curvature at $v$ is the expected genus of a random surface in $S(v)$. }

\paragraph{}
The theorem generalizes from manifolds {\bf Dehn-Sommerville q-manifolds}. These
are $q$-dimensional finite abstract complexes for which all $S(x)$ are 
Dehn-Sommerville $(q-1)$-manifolds of Euler characteristic $1-(-1)^q$. 
{\bf Euler's gem } implies that $q$-manifolds are q-Dehn-Sommerville manifolds. A
{\bf Dehn-Sommerville q-sphere} is a Dehn-Sommerville q-manifold of Euler characteristic  $1+(-1)^q$.
The sphere formula $\sum_{x \in G} \omega(x) \chi(S(x))=0$ implies that every
odd-dimensional Dehn-Sommerville manifold is a Dehn-Sommerville sphere.
\index{Dehn-Sommerville manifold}
\index{Dehn-Sommerville sphere}
\index{Euler's gem}

\satz{A level subset defined by a function $g: V(G) \to \{0,\dots, k\}$
on a Dehn-Sommerville $q$-manifold $G$ is a $(q-k)$ Dehn-Sommerville manifold or empty. }

\paragraph{}
It implies :

\satz{For odd-dimensional Dehn-Sommerville complexes, $j_g(v)=0$ for all $v$.
For even dimensional Dehn-Sommerville complexex $j_g(v) = (1-\chi(S_g(x))/2$. }

\paragraph{}
Let $(\Omega,\mathcal{A},P)$ be a probability measure on locally injective functions 
with the property that $Tg=-g$ is an automorphism of the probability space. We call
the index expectation curvature $k(v)={\rm E}[j_g(v)]$ a {\bf symmetric curvature}.
The Levitt curvature and the color curvature are both symmetric. We have in particular:

\satz{Odd-dimensional Dehn-Sommerville manifolds are flat.}

 \vfill \pagebreak
\begin{center}\scalebox{0.8}{\includegraphics{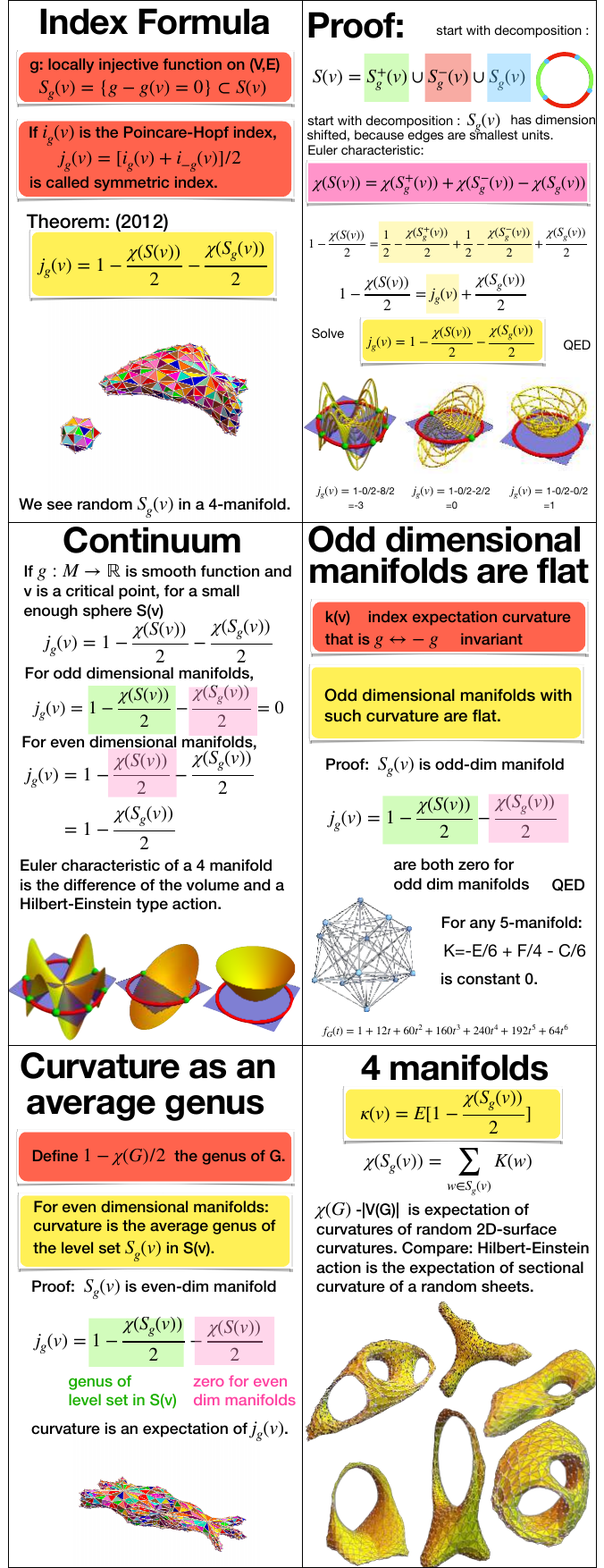}}\end{center} \pagebreak
\setcounter{section}{11} \chapter{Unit 11: Quadratic cohomology }

\summary{
If $G$ is a finite abstract simplicial complex, the quadratic characteristic 
$\chi_2(G) = \sum_{x \cap y \in G} \omega(x) \omega(y)$ is a quadratic version 
of Euler characteristic $\chi(G) = \sum_x \omega(x)$. It comes with 
a cohomology similarly as $\chi(G)$ comes with simplicial cohomology. For
$q$-manifolds with boundary $\chi_2(G)=\chi(G)-\chi(\delta G)$. }

\paragraph{}
Classical theorems like Gauss-Bonnet, Poincare-Hopf, Euler-Poincare or Brouwer-Lefschetz are 
related to Euler characteristic $\chi(G) = \sum_{x \in G}  \omega(x)$, both in the continuum as 
in the discrete. There are higher order theories
for which {\bf quadratic characteristic} $\chi_2(G) = \sum_{x \cap y \in G} \omega(x) \omega(y)$ is the 
simplest. Some theorems generalize directly but there are changes, the most important being that 
$\chi_2$ is of a topological nature. It is not homotopy invariant. If $G$ is the closure of $\{x\}$ 
then $\chi_2(G)=\omega(x)=(-1)^q$. For manifolds without boundary, $\chi_2(G)=\chi(G)$. 
\index{quadratic characteristic}
\index{Wu characteristic}
\index{quadratic cohomology}

\paragraph{}
{\bf Quadratic cohomology} works with functions on {\bf pairs of intersecting simplices}
$(x,y) \in G \times G$. The {\bf exterior derivative} is defined as
$df(x,y) = d_x f(x,y) + \omega(x) d_y(f,y)$, where $d_x,d_y$ are the partial
simplicial exterior derivatives.  
Cubic, quartic and higher characteristics and cohomology are defined similarly. 

\paragraph{}
With ${\rm dim}(x,y)={\rm dim}(x)+{\rm dim}(y)$, let
$\Lambda_{k}= \{ f: G_k \to \mathbb{R}$, where $G_k=\{ (x,y), {\rm dim}(x,y)=k \}$ is
the space of $k$-forms and let $f_k$ denote its dimension. Then $w(G) = \sum_{k=0}^{2q} (-1)^k f_k$. The space
$\Lambda = \oplus_{k=0}^{2q} \Lambda_k$ of {\bf differential forms} has as dimension the number of 
ordered pairs $(x,y) \in G^2$ with $x \cap y \in G$. It is a finite dimensional vector 
space and ${\rm str}(1) = \chi_2(G)$ is the quadratic characteristic. The cardinalities of the 
quadratic complex is encoded by 
$f_{kl} = \{ (x,y), x \cap y \neq \emptyset, {\rm dim}(x)=k, {\rm dim}(y)=l \}$. 

\paragraph{}
The Dirac operator $D=d+d^*$ and the Hodge operator $H=D^2=d d^* + d^* d$ are defined as before.
If $G$ has maximal dimension $q$, then $H$ has $2q+1$ blocks $H_k$. It belongs to forms
$f(x,y)$ acting on pairs $(x,y)$ of simplices with ${\rm dim}(x) + {\rm dim}(y) = k$. 
The kernel of $H_k$ is called the $k$'th Betti number $b_k$.
As in the linear case, we can split $\Lambda(G)$ into {\bf even forms} $\Lambda_{2j}$ and
{\bf odd forms} $\Lambda_{2j+1}$. The matrix $D$ is an isomorphism between even and odd forms 
on the image of $H$ so that ${\rm str}(H^n)=0$ for $n>0$. The analog of Euler-Poincar\'e is
proven in exactly the same way using heat deformation.
\index{quadratic Dirac operator}
\index{quadratic Hodge operator}

\satz{$\chi_2(G) = \sum_{k=0}^{2q} (-1)^k f_k = \sum_{k=0}^{2q} (-1)^k  b_k$.}

\paragraph{}
An automorphism $T$ of $G$ also induces an automorphism on the intersection diagonal of 
$G \times G$ by $T(x,y) = (T(x),T(y))$ the reason being that if $x \cap y \neq \emptyset$,
then also $T(x) \cap T(y) = T(x \cap y) \neq \emptyset$. The Lefschetz number $\chi_{2,T}(G)$
is as in the linear case the super trace of the Koopman operator $U$ on cohomology. The 
index of a fixed point $(x,y)$ is $i_T(x,y) = \omega(x) \omega(y) {\rm sign}(T|x) {\rm sign}(T|y)$.  
Also the Lefschetz fixed point theorem generalizes with the same heat proof to a 
{\bf quadratic Lefschetz fixed point theorem}:

\satz{$\chi_{2,T}(G) = \sum_{(x,y) \in F} i_T(x,y)$.}
\index{quadratic Lefschetz fixed point theorem}

\paragraph{} 
It generalizes the Euler-Poincar\'e theorem because for $T=Id$, $\chi_{2,T}(G) = \chi_2(G)$
and every $(x,y)$ with $x \cap y \neq \emptyset$ is then a fixed point of $T$ with 
index $\omega(x) \omega(y)$ so that the right hand side in the Lefschetz theorem is
the definition of $\chi_2(G)$. 

\paragraph{}
Also Gauss-Bonnet, Poincar\'e-Hopf and index expectation generalize.
Curvature $K(v)$ for quadratic valuation is
$K(v) = \sum_{x \sim y} \omega(x) \omega(y)/|x|$.

\satz{$\chi_{2}(G)=\sum_{v \in V} K(v)$.}

\paragraph{}
For manifolds with boundary, we know how to compute it.

\satz{ $\chi_2(G) = \chi(G) - \chi(\delta G)$.  }

For even $q$, this gives $\chi_2(G)=\chi(G)$ even with boundary.
For odd $q$, this gives $\chi_2(G)=-\chi(\delta G)/2$.
\index{boundary theorem for quadratic cohomology}

\paragraph{}
Examples: 
1) For $P_2=\{ \{1\},\{2\},\{3\},\{1,2\},\{2,3\} \}$, there 15 interacting pairs $(x,y)$ where
$f_0=3,f_1=8,f_2=4$. The Wu characteristic is $3-8+4=-1$.
2) The utility graph $G=K_{3,3}$ has Euler characteristic $\chi(G)=-3$ and Wu characteristic  
$\chi_2(G)=15$. It is maximal among all graphs with 6 vertices. 
3)  For a star graph with $n$ rays, the Wu characteristic is $n^2-3n+1$. For example, 
for $n=0$, it is $1$, for $n=1$ it is $-1$ for $n=9$ it is $55$. All star graphs are
contractible illustrating the non-homotopy invariance. 
4)  Adding a hair to a 2 sphere reduces the Wu characteristic by $2$.  
5) The Wu characteristic of the cube graph is $20$, the Wu characteristic of
the dodecahedron is $50$. Both graphs have no triangles and constant vertex degree $d=3$
so that in both cases, the curvature is constant $5/2$. 
6) The dunce hat has quadratic cohomology $b=(0,0,0,1,2)$. 
7) The Moebius strip has $b=(0,0,0,0,0)$ while the cylinder has $b=(0,0,1,1,0)$. The cohomology 
sees the twist. The Lefschetz fixed point theorem implies that if an automorphism preserves 
harmonic 2-forms and reverses harmonic 3-forms, then there are at least 2 fixed points. We
once speculated that the preservation of harmonic 2-forms could be related to area preservation
and the twisting of harmonic 3 forms be related to twist so that quadratic cohomology could
give Birkhoff's fixed point theorem.

 \vfill \pagebreak
\begin{center}\scalebox{0.8}{\includegraphics{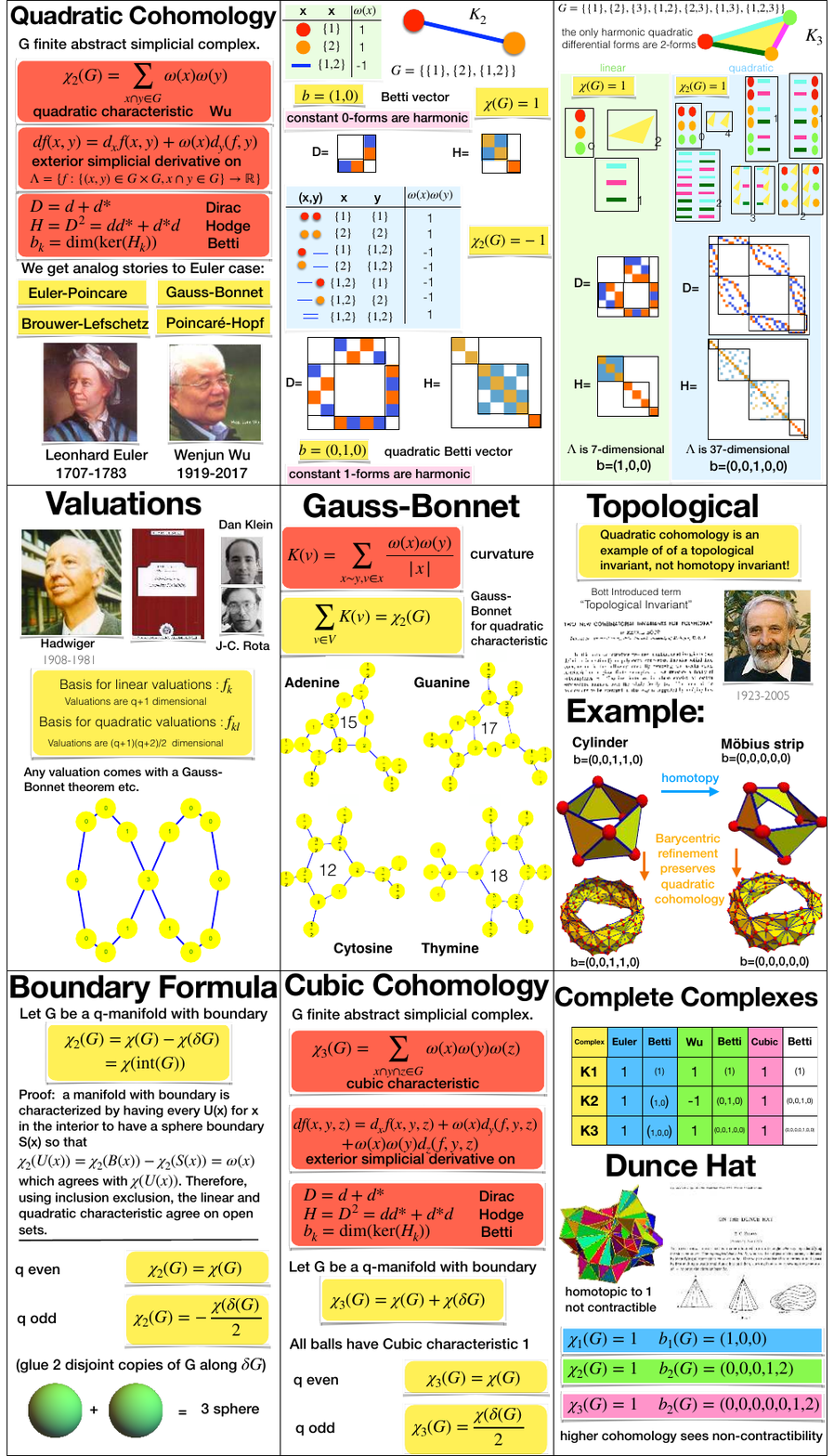}}\end{center} \pagebreak
\setcounter{section}{12} \chapter{Unit 12: Higher Green }

\summary{
The $k$-point {\bf Green function} 
$g_m(x_1, \dots, x_k)=\prod_{j=1}^k \omega(x_j) \chi_m(\bigcap_{j=1}^k U(x_j))$
for the {\bf $m$'th characteristic} 
$\chi_m(A) = \sum_{x \in A^m, \bigcap_j x_j \in A} \prod_{j=1}^m \omega(x_j)$ of $A \subset G$ 
satisfies $\chi_m(G) = \sum_{x \in G^k} g_m(x_1, \dots, x_k)$.  }

\paragraph{}
The {\bf $m$'th characteristic} is defined as
$\chi_m(G) = \sum_{{\bf x} \in G^m, \bigcap {\bf x} \in G} \omega(x)$.
It is the Euler characteristic in the linear case $m=1$ and the quadratic characteristic 
in the quadratic case $m=2$. If used vector notation ${\bf x}=(x_1, \dots, x_m)$
and multi-index notation $\omega({\bf x}) = \prod_{j=1}^m \omega(x_j)$ and 
$\bigcap {\bf x} = x_1 \cap x_2 \cap \cdots \cap x_m$. 
\index{higher characteristic}
\index{multi-index notation}

\paragraph{}
The {\bf $k$-point Green function} on $G^k$ is $g_m({\bf x})=\omega({\bf x}) \chi_m(\bigcap U({\bf x}))$,
where $U(x_j)$ is the star of $x_j$ and the star $U({\bf x})$ of ${\bf x} \in G^k$ is the intersection
of the stars of $x_1, \dots, x_k$. The $k$-tuple ${\bf x}=(x_1, \dots, x_k)$ does not need
to intersect but the $k$-point Green function are still local because $U({\bf x})$ is empty, if one of
the ``particles" $x_j$ is out of reach. 
\index{k-point Green function}
\index{Green function}

\paragraph{}
The Green function identity $\chi(G)=\sum_{x,y} g(x,y)$ for the quadratic case $m=2$ generalizes
but already in the cubic case $g(x,y,z)$ is a tensor. For $m=2$, we could think of $g(x,y)$
as a matrix, where it was the inverse of $L(x,y)$. 
While $\chi_m(G)$ involves only $m$-tuples of $G$ that intersect, the Green functions 
$g(x_1, \dots, x_k)$ involve all $k$-tuples of $G$. Think of a {\bf $k$-point potential energy}.
The energy theorem generalizes to

\satz{ $\chi_m(G) = \sum_{{\bf x} \in G^k} g_m(\bf{x})$.}
\index{energy theorem}

\begin{proof} 
Generalize to energized complexes where $\omega({\bf x})$ is replaced with 
a general function $h_m(x_1,\dots, x_m)$ of $m$ variables producing a m-linear
valuation.  Since $h_m \to g_k$ is linear, we only need to verify the
statement in the simplest possible case, where $h$ is $1$ only for a single
configuration ${\bf z}=(z_1, \dots, z_m)$ and $0$ else. The left hand side is then $1$. 
On the right hand side, we have to look at all 
${\bf x}=(x_1,\dots, x_k)$ for which ${\bf z} \in U({\bf x})=U(x_1) \cap U(x_2) \cap \cdots \cap U(x_m)$.
$\sum_{{\bf x} \in G^k} \omega({\bf x}) \chi_m(U({\bf x})) = \prod_{j=1}^k \sum_{x_j \in G} \omega(x_j) 
             \chi_m(U(x_j)) =1$.
\end{proof}

\paragraph{}
A second major point is the {\bf sphere formula} for 
$S({\bf x}) = \delta U({\bf x})=B(\bf x) \setminus U({\bf x})$,
where we write similarly as $U({\bf x})=\bigcap_j U(x_j)$ also $B({\bf x})=\bigcap_j B(x_j)$.  
\index{sphere formula for higher characteristic} 

\satz{$\sum_{{\bf x} \in G^k} \omega({\bf x}) \chi_m(S({\bf x})) =0$.}

\paragraph{}
This means that the total m'th characteristic of all unit spheres of even $k$-point configurations
is the same than the total m'th characteristic of all unit spheres of odd $k$-point configurations.

\begin{proof}
The energy theorem also works of $U({\bf x})$ is replaced by $B({\bf x})$, which is the
closure of $U({\bf x})$. The two equations
$0 = \sum_{{\bf x} \in G^k} \omega({\bf x}) \chi_m(U({\bf x}))$ and
   $0 = \sum_{{\bf x} \in G^k} \omega({\bf x}) \chi_m(B({\bf x}))$
and the local valuation formula blow prove the theorem.
\end{proof}

\paragraph{}
For characteristic $m=1$, the formula $\chi( \overline{U}) = \chi(U) + \chi(\delta U)$ follows
from the fact that every simplex in $\overline{U}$ is either in $U$ or $\delta U$. The 
{\bf local valuation formula} is:

\satz{ $ \chi_m(B({\bf x})) = \chi_m(U({\bf x})) - (-1)^m \chi_m(S({\bf x}))$. }
\index{local valuation formula}

\paragraph{}
For odd $m$ this means  $\chi(U) + \chi(\delta U) = \chi( \overline{U})$. For
quadratic characteristic $m=2$, one has $\chi_2(U(x)) - \chi_2(S(x)) = \chi_2(B(x))$. 
If $B(x)$ is a replaced with a manifold $G$ with boundary and we chose $k=1$, one gets the boundary formula
$\chi_2(G) = \chi_2({\rm int}(G)) - \chi_2(\delta G)$ which then simplifies to $\chi_2(G) = \chi(G)-\chi(\delta G)$. 

\satz{ If $G$ is manifold we have $\chi_m(G) = \chi(G)-(-1)^m \chi(\delta G)$. }
\begin{proof}
Use induction and that the left hand side is $\chi_m({\rm int}(G)) - (-1)^m \chi_m(S(\delta G))$.
\end{proof} 

\paragraph{}
Higher characteristics of manifolds with boundary can be expressed using Euler characteristic.
For even-dimensional manifolds with boundary $\chi_m(G) = \chi(G)$. For odd-dimensional manifolds with
boundary $\chi_m(G) = -(-1)^m \chi(\delta G)/2$. For even dimensional balls $G$ for example, where the boundary 
is an odd-dimensional sphere, $\chi_m(G)=1$ independent of $m$, while for odd-dimensional balls $\chi_m(G)=(-1)^{m+1}$.  

\paragraph{}
Higher characteristic $\chi_m(G)$ can be written as $\chi_m(G,G,\dots, G)$ if one defines the multi-linear
expressions $\chi_m(A_1, \dots, A_m) = \sum_{x_1 \cap \dots \cap x_m \in G, x_i \in A_i} \omega({\bf x})$.
Like for multi-linear forms in linear algebra, fixing all but one coordinates produces a valuation. 
Like Euler characteristic had the $f$-vector and quadratic valuations had the $f$-matrix $f_{kl}$, 
multi-linear forms define the {\bf $f$-tensor} $f_{k_1, \dots, k_m}$. Applied to $m$ Dehn-Sommerville 
vectors, we get Dehn-Sommerville invariants for manifolds answering a question of Gruenbaum. 
The most prominent Dehn-Sommerville invariant is $X=(1,-1,1,-1, \dots )$ which produces characteristic 
invariants $f(X,X,\dots, X)$. Let us abbreviate $m$-valuation for multi-linear valuation of $m$ variables:

\satz{ $\chi_m$ is the only $m$-valuation that is Barycentric invariant.}
\begin{proof}
$(1,-1,1,-1, \dots, \pm 1)$ is the only eigenvector of the {\bf Barycentric operator}.
\end{proof} 

Quantities with this property are {\bf topological invariants}. 
\index{topological invariant}
 \vfill \pagebreak
\begin{center}\scalebox{0.8}{\includegraphics{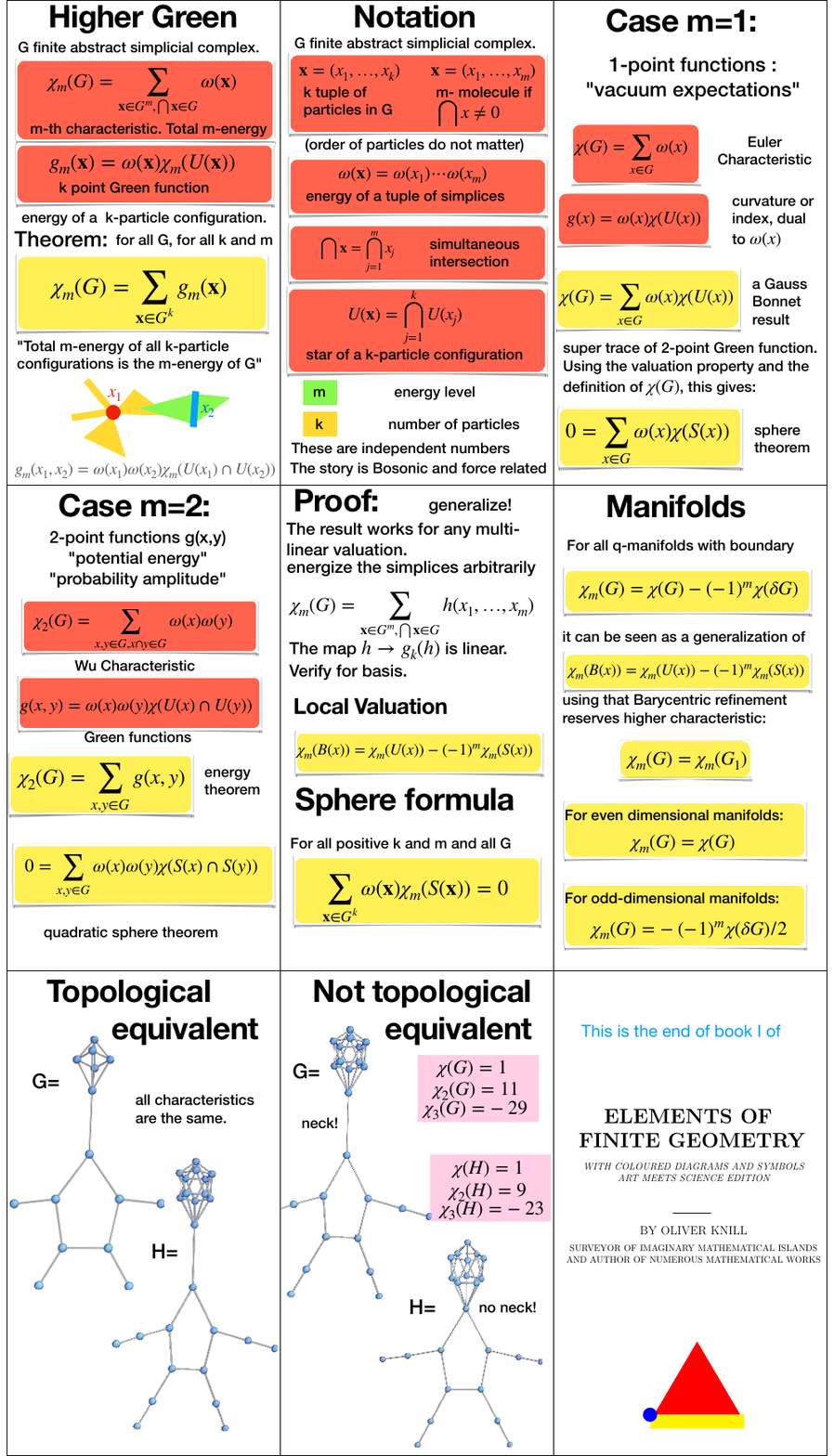}}\end{center} \pagebreak
                         \section*{Footnotes}

This section contains some remarks to the text. 
Brevity in the main text was not to get bogged down. There were previous attempts to summarize 
things and they always hit a wall. I gave myself a self-imposed limit of 2 pages per unit and a one 
minute limit for the accompanying videos.  The illustrations are from slides written for such youtube shorts
written during the summer. 
Texts were drafted partly during the earlier summer during travels. It started in the Swiss alps 
with a view onto the Mischabel mountain range seen on the title page. One can see about 12 mountain peaks.
The footnotes and references aim less at optimizing 
length and therefore meander more. Possible in the future could be an appendix with examples and exercises. 
But it will be more important to work on the next batch of 12 units, if circumstances permit. This summer
was the first in 24 years, where I had not been teaching in summer school, providing an opportunity to
work on this project.  \\

{\bf 0. Introduction} \\

The motivation on why to focus on finite structures is a reaction to
the incompleteness theorem \cite{Goedel1940} which stresses that the foundations of mathematics will
always remain uncertain if the assumptions are strong enough. 
Popularized masterfully in \cite{Hofstadter}, it had been one of the most pivotal high school experiences that 
drew me personally to mathematics. I got exposed to ``strict finitism" by Erwin Engeler at ETH 
(he taught us a set theoretical "topology" course and "mathematical software" course). 
He mentioned strict finitism once in person but \cite{EnglerStrictFinitism} is one 
of Engeler's articles from 1971.
Ernst Specker (from whom I took two semesters of linear algebra, as well as
logic and model theory courses and also participated in various seminars), mused about the relation between finite
mathematics and the continuum and informed the class once about the work of Paul Kustaanheimo 
\cite{Kustaanheimo}, an attempt to do mathematics in finite fields, rather than the real numbers. Kustaanheimo
collaborated with Eduard Stiefel while visiting Zuerich in 1964, certainly knowing Specker, as Specker and Stiefel
both were students of Heinz Hopf. Peter L\"auchli, a student of Specker (from whom I took two semesters of 
"calculus" (both single and multi-variable) and a logic course about forcing) gave later a mesmerizing course on 
non-standard analysis based on Nelson's internal set theory \cite{Nelson77,Nelson}. See
\cite{GoldenRotations} for an other personal story about Specker and L\"auchli (and other ETH teachers).
The textbook that L\"auchli recommended was \cite{Robert} and he followed it partly.
IST an effective link between finitist and traditional mathematics
as it just changes language and not structures. It shows that compact spaces can be treated like finite spaces.
Nelson also reformulated probability theory \cite{Nelson} and called it``radically elementary". 
Finite models were also pushed by other of my teachers. J\"urgen Moser mentioned that finite systems show similar
KAM features than the continuum \cite{Ran74} and my PhD advisor
Oscar Lanford II worked in his later career on finite approximations of dynamical systems \cite{lanford98}.
Finally, finitism is a honest reaction to the fact that all 
data we ever acquire and process are finite. It is also a computer science approach to mathematics. 
The Taylor formula mentioned in the introduction appeared in a Pecha-Kucha talk \cite{ArchimedesFunctions},
a presentation in which 20 slides are shown with a strict $20$ second limit to talk for each.
The discrete Taylor theorem is attributed to James Gregory, 
a contemporary of Newton. One can find discrete calculus also in the original writing of Leibniz.
I used discrete calculus excessively in single variable Math 1a courses from 2011 to 2024. I was
fortunate to have been given the freedom to experiment with a fresh take on single variable calculus then. 
\cite{VanderWaerdenAlgebra} contains the formula on page 69 in the context of interpolation. 
The topics discussed in these notes (and beyond) have previously seen attempts to be summarized:
\cite{knillcalculus,KnillILAS,KnillBaltimore,AmazingWorld,LusternikSchnirelmann2024,DehnSommervilleManifolds,Simon2026}. \\
\index{Goedel}
\index{incompleteness theorem}
\index{finitism}
\index{strict finitism}
\index{Ernst Specker}
\index{Peter L\"auchli}
\index{Erwin Engeler}

{\bf 1. Gauss-Bonnet} \\

For the discrete Gauss-Bonnet theorem, see \cite{Eberhard1891,Levitt1992,elemente11,dehnsommervillegaussbonnet,
cherngaussbonnet,valuation,DehnSommerville}. Norman Levitt, was the first who has introduced curvature of complexes
in higher dimensions. He did not stress it as a Gauss-Bonnet theorem. The integral geometric interpretation as
an index expectation verifies with Nash's embedding theorem that the Gauss-Bonnet-Chern curvature for even 
dimensional manifolds is a limiting case of this discrete measure. Many other curvatures have been proposed
graph theory, some using probability theory \cite{Ollivier,JostLiu}. For me, a curvature is interesting if it 
satisfies Gauss-Bonnet. Also notions like sectional curvature should 
be curvatures of a two dimensional geodesic surface passing through a point at which the curvature is measured.
I learned about Levitt's paper only years after \cite{cherngaussbonnet}. I had seen the general formula only after
painfully going through small dimensional manifold cases, where Dehn-Sommerville relations obscure
a bit the patterns. This already happens in the case of oldest formula $K=1-{\rm deg}(v)/6$ for triangulated
two dimensional surfaces, where the Dehn-Sommerville identity $2|E|=3|F|$ hides that $K=1-{\rm deg}(v)/2+{\rm deg}(v)/3$. 
This curvature formula appears in \cite{princetonguide} after writing \cite{elemente11}.
A crucial problem had been for then to explain why we always see constant zero curvature 
for odd-dimensional manifolds. This turned out to be a rich problem. It motivated 
to look at integral geometric approaches and a produced several encounters with the Dehn-Sommerville story.
The integral geometric connection proves that the curvature is the real deal. 
The discrete curvature is not just some ``analogy".  It produces the Gauss-Bonnet-Chern integrand if we 
take a discrete approximation of a compact Riemannian manifold. 
Related more general curvatures appeared in \cite{Higuchi,FormCurvatures}.
For the classical case, see \cite{Wu2007}. The classical Gauss-Bonnet-Chern evolved over decades:
\cite{HopfCurvaturaIntegra,Allendoerfer,Fenchel,AllendoerferWeil,Chern44,chern1990,Cycon}
Discretisations of space have appeared long before Poincar\'e, but less for general manifolds. Historically,
one has focussed more on polytopes \cite{Ziegler,gruenbaum}. ``Monsters" as Lakatos described counter examples
to the Gem formula have made geometers more careful.  It is a ``delicate task" \cite{Devadoss}.
The result was that discussions about polyhedra are
mostly limited to convex objects. If we discard the continuum, we need to be extra clear, in order
not to repeat the scandal which Lakatos described masterfully. If we have a finite combinatorial object without
any embedding into some Euclidean space, we have to be precise with what is considered a ``face" of a given dimension.
For history see \cite{Dieudonne1989}.  For discrete geometries see \cite{GSD,Regge,MarsdenDesbrun} for computational geometry 
\cite{Devadoss,ComputationalGeometry}, for integrable systems \cite{BobenkoSuris}, 
for computer graphics or computational-numerical methods \cite{Bobenko,Bobenko2008,CompElectro2002}. 
\index{Gauss-Bonnet-Chern}
For graph theory in general see \cite{Diestel,HararyGraphTheory,Merris,BR}
For combinatorics of simplicial complexes in general \cite{May1967,Stanley86,RhodesSilva}.

{\bf 2. Poincar\'e-Hopf} \\

The discrete Poincar\'e-Hopf appeared first in \cite{poincarehopf}. It was useful for
\cite{KnillEnergy2020}.  The work with Frank Josellis in \cite{josellisknill} led to simplifications. We
hope to get to Morse and Lusternik-Schnirelmann topics in the second volume.
The parametrized version appeared in \cite{parametrizedpoincarehopf}.
More elaborations came in \cite{MorePoincareHopf} and \cite{PoincareHopfVectorFields}. 
A poetic approach to numbers is obtained by looking at numbers geometrically \cite{PrimesGraphsCohomology}.
I had then not been aware of \cite{Bjoerner2011} which introduced that complex already.
\cite{PrimesGraphsCohomology} however allows to interpret the Riemann hypothesis (the Mertens connection)
from a topological point of view. For an early discrete version in two dimensions based on an embedding
in a surfaces, see \cite{Glass1973}. I learned about Glass'theorem from Kate Perkins thesis
at Harvey Mudd \cite{KatePerkins}.
The speed with which the simplex generating function can be computed recursively using the formula
given in the text, needs to be studied more. Computing the number of simplices in a graph is 
known to be a hard problem. The clique problem is NP complete \cite{Karp}. What the recursive formula does is 
relegate the problem to subgraphs in unit spheres.
To the classical story of Poincar\'e-Hopf: 
Poincar\'e proved the index theorem in chapter VIII of \cite{poincare85}. Hopf extended it to
arbitrary dimensions in \cite{hopf26}. For manifolds with boundary, see \cite{Pugh1968}. 
\cite{Spivak5} contains a modern exposition.
The collection \cite{HistoryTopology} contains in particular an article on the history of 
graph theory. A story about Euler characteristic and polyhedra is told in \cite{Richeson,coxeter}.
See \cite{Hatcher,Spanier,Rotman} for
notations in algebraic topology, \cite{HararyGraphTheory,Biggs,BM} for graph theory.
Valuations have been pioneered by Klee \cite{Klee63} and Rota \cite{Rota71} who already saw that
there is a unique normalized valuation. As for Hadwiger see \cite{hadwiger, Klain1995}.
\index{Kate Perkins}
\index{Frank Josellis}
\index{Glass theorem}
\index{Mertens conjecture}
\index{Morse theory}
\index{Lusternik-Schnirelmann theory}

{\bf 3. Index expectation} \\

The Buffon needle problem, stated and solved in 1777, is one of the
first problems in the theory of geometric probability. The subject is also called integral geometry.
$L = \pi {\rm E}[n]$ gives the length if lines are distance $2$ apart. 
See \cite{Santalo1}.  Crofton's book on local probability is \cite{Crofton1968}.
More general references on integral geometry are \cite{Santalo,KlainRota,Schneider1,Blaschke,Nicoalescu}.
The first time, it seems that curvature was seen as an expectation of indices is in the work of Banchoff
\cite{Banchoff1967,Banchoff1970}.
It has also appeared in graph theory \cite{IntegralGeometryGraphs}.
As for index expectation, see \cite{indexexpectation} 
\cite{indexformula},     
\cite{ConstantExpectationCurvature} 
\cite{DiscreteHopf}    
\cite{DiscreteHopf2}.   
For the Nash embedding theorem source, see \cite{EssentialNash}. Gauss-Bonnet-Chern is older than
Chern, if one allows embeddings.  The first intrinsic proof of GBC without embedding the manifold
is \cite{Chern44}. For historical remarks see \cite{chern1990}. I myself learned GBC first from 
Patodi's proof \cite{Cycon}.  It was Konrad Osterwalder who made me look into chapter 12 of the book of 
Cycon-Froese-Kirsch and Simon to assist in a seminar on Patodi’s proof of the Gauss-Bonnet-Chern theorem. 
I spent a week alone on the alp “Salmenfee” in the Swiss mountains during a summer break, working through 
that chapter. The seminar eventually did not take place; the exposure to that book however turned 
out to be life changing for me, as the book contained a nice concise introduction into modern Riemannian geometry
and interesting spectral theory. Integral geometry is extremely useful as it links the discrete with the continuum. 
Because Poincar\'e-Hopf indices are divisors, there is a direct link between the 
discrete and continuum. Integral geometry allows to bypass any discrete Riemannian
manifold set-ups. Distances can be implemented using embeddings for example. 
Index expectation was crucial for the product formula for the Shannon product
\cite{Shannon1956}. See also \cite{ComplexesGraphsProductsShannonCapacity}. \\
\index{Buffon needle problem}
\index{Crofton's formula}

{\bf 4. Euler's Gem} \\

The task to define what a manifold is combinatorially essentially boils down to the task 
to define ``spheres", a question going back to Herman Weyl \cite{Weyl1925}, at a time
when the ``Grundlagen crisis" had been raging still. (There is still a Grundlagen crisis today 
of course but similarly as the possibilities of a nuclear armagedon danger is blended out from
the public consciousness, mathematicians usually ignore the possibility of an inconsistency to 
occur. They have doubled down even by using structures that go beyond ZFC like Grothendieck 
universes to avoid objects like the category of all categories. 
There are two important early approaches to discrete manifolds. For simplicial complexes, see
\cite{JonssonSimplicial, Stanley86,Stanley1996}. 
One approach is by Robin Forman \cite{forman95,forman98,Forman1999,Forman2003}
who built ``discrete Morse theory", the other by Ivashchenko (Evako) \cite{Ivashchenko1993,I94,CYY}
whose work is sometimes placed into a field called ``digital topology". Forman defines spheres 
Morse theoretically using functions, Evako uses induction. The compact version done here has 
been developed by us, first in the graph theoretical set-up, later in abstract simplicial 
complexes using the Alexandroff topology \cite{Alexandroff1937,May2008,KnillTopology2023}. 
What is important to make the definition decidable
is to use "contractibility" (sometimes called "collapsible") and not "homotopic to 1" in the
definition of a sphere. With the later, one would run into computer science trouble: there would
be no Turing machine which could decide whether a given structure is a manifold or not. 
With "contractibility" as defined, there is a definite algorithm which terminates polynomial time
once the dimension is fixed and decides whether a given graph or a given finite abstract simplicial
complex is a manifold or not. One of the simplest complexes that is homotopic to 1 but not contractible
is the dunce hat \cite{ZeemanDunceHat}.
Discrete Reeb \cite{knillreeb} shows equivalence of Morse and Recursive approaches.
A good general reference for the topic of Euler's gem is \cite{Richeson}.
I gave a mathtable talk about this on 2/6/2018 to honor Euler's day 2/7/18 \cite{KnillEulerGem}
The innocent looking Euler Gem formula is infamous for having been proven
wrong repetitively, and not by amateur mathematicians or but by the best mathematicians in the world. 
It so has become a show case for the epistemology of mathematical discovery \cite{lakatos}. 
I myself disagree with Lakatos that science is an evolutionary approach to truth. It should
be true from the beginning. The culprit usually is that the definitions are
not good. The algebra of the join operation goes back to \cite{Zykov} in graph theory. It might even
have appeared already in the Principia Mathematica from 1910. 
The analog for simplicial complexes is to take $A \cup B \cup \{ a \cup b, a \in A, b \in B \}$ for
the join of two closed sets.  It
produces a monoid structure on complexes and spheres are a sub-monoid. Also Dehn-Sommerville spaces
form a submonoid.  The classification of regular polytopes is usually done in the continuum in the context of
polytopes \cite{cromwell,symmetries,Sutton,Ziegler,gruenbaum}, which invoke Euclidean space and convexity. \cite{Gruenbaum2003} discussed
the difficulty with the term ``polytop". 
For \cite{coxeter},a polytop is defined as as a convex body with polygonal faces. 
For source work, \cite{Schlafli,Schoute,AliciaBooleStott}. The origin of Euler characteristic is captivating
\cite{Aczel}. Polytop definitions are given in \cite{Schlafli,coxeter,gruenbaum,symmetries}. 
Topologists started with new definitions \cite{alexandroff,Fomenko,ConwayCapstone,Spanier}.
Abstract simplicial complexes appeared in 1907 \cite{DehnHeegaard}. See
\cite{BurdeZieschang,MunkholmMunkholm}. \cite{alexandroff} uses the name 
unrestricted skeleton complex. In \cite{Whitehead} there are called
symbolic complexes. For references on simplicial complexes, see 
\cite{JonssonSimplicial, Stanley86,Stanley1996} but definitions can vary. Especially the treatment 
of $\emptyset$. Insisting on having simplicial complexes not contain $\emptyset$ is much, much 
more natural and makes the theory look like the theory we all know when studying differential 
topology or differential geometry.  \\
\index{Evako, Ivashchenko}
\index{Robin Forman}
\index{Discrete Morse theory}

{\bf 5. Euler-Poincare} \\

Finite exterior calculus has been used already by the founders of algebraic
topology. The one-dimensional discrete Laplacian $H_0$ is now called the {\bf Kirchhoff matrix} \cite{Kirc}.
Discrete Hodge theory appeared first in \cite{Eckmann1944}. One of the talks
given by Beno Eckmann at ETH given at the University of Z\"urich that focussed on Euler characteristic 
made me personlly curious about this topic.  Looking at higher order Laplacians is pretty established 
also in the discrete (an example \cite{DanijelaJost}). Discrete calculus has been re-invented again and again
especially in the applied and computer science literature. Examples are \cite{marsden_discretecalculus,marsden_poincarelemma}.
Discrete McKean-Singer \cite{McKeanSinger} appeared in \cite{knillmckeansinger}. Taking simple ideas going
back to Dirac and Hodge (they were office neighbors who did not talk much together as
Michael Atiyah once happily pointed out in an interview as this resulted him building a career combining both their
ideas).  Discrete calculus has been reinvented again and
again, in later time, it has appeared mostly in applied or computer science literature. 
\cite{Edelsbrunner}. 
The main definition of using linear algebra in topology have already been forged 
by Poincar\'e or Betti. The matrix $D$ is usually given in terms of its blocks, which are called 
{\bf incidence matrices}. It would have been so much easier already at the time of Poincar\'e to
combine these matrices to a common matrix $d$ on all forms and to define $D=d+d^*$ and then study
the matrix $D^2$. See \cite{KnillILAS}. The insight that one does not have to look at 
anti-symmetric functions, but just can look at functions, once an orientation 
is fixed on simplices, is probably due to Hassler Whitney who used this especially in the context 
of the cup product. To bypass the chain complex notations and use linear algebra is due to Hodge. 
It is such a relieve especially when working with computer algebra implementations not to have to deal
with equivalence classes, but simply with null-spaces of matrices.
The concept of the square root of the Laplacian is due to Dirac. Clifford algebra constructions 
show that we can bypass awkward Dirac matrix constructions. 
The Hodge definition for cohomology is more intuitive than equivalence classes 
"closed forms modulo exact forms". It is a good (a bit challenging but doable) exercise for 
an introduction linear algebra course to show that the two pictures equivalent. 
About the notation, we use $H$ for Hodge and because we want to use $L$ for the 
connection Laplacian. We write {\bf Dirac operator} and not {\bf Dirac matrix} in order
not to confuse with Dirac matrices in mathematical physics, which are used to get the 
square root of the Laplacian. \\
\index{Laplacian}
\index{Kirchhoffmatrix}

{\bf 6.  Unimodularity} \\

The result was discovered in December 2015 while experimenting 
with quadratic cohomology. Neither any discrete nor continuum analog of this story
seems have appeared before. 
It was presented on a mathtable talk "Bowen-Lanford Zeta function"
on 10/18/2016 \cite{KnillBowenLanford}, referring to \cite{BowenLanford}.
The proof appeared in \cite{Unimodularity}. Some extension appear in \cite{CountingMatrix}
like when $\omega(x)$ is replaced with $1$.
A reason for the delay was a summer excursion into number theory \cite{Experiments}.
The result still surprises me today. A more direct alternative proof of unimodularity 
was given in \cite{MukherjeeBera2018}. I still feel the unimodularity result is remarkable and
I'm glad I spent some effort to have it published \cite{KnillEnergy2020}. Polishing something for
an actual publication uses a lot of time and often is a fool's errand anyway, if one deals with 
a topic, that has not been covered before. There were various generalizations
\cite{GreenFunctionsEnergized},
\cite{EnergizedSimplicialComplexes},
\cite{EnergizedSimplicialComplexes2}, and
\cite{EnergizedSimplicialComplexes3}. For some spectral problems in the 
context of Dirac and connection matrices, see \cite{ConnectionDirac}. \\
The connection story is related to the {\bf Szpilrajn-Marczewski theorem}: any finite simple graph 
$A$ can be realized as a connection graph of a finite set $G$ of non-empty sets 
\cite{Szipilrajn-Marczewski,ErdoesGoodmanPosa} as pointed out in \cite{ComplexesGraphsProductsShannonCapacity}. \\
\index{Szpilrajn-Marczewski theorem}

{\bf 7. Lefschetz fixed point formula} \\

The discrete Brouwer Lefschetz theorem appeared in \cite{brouwergraph}.
My proof there was close to Hopf's approach \cite{Hopf28} which also is in
\cite{alexandroff} (Alexandrov and Hopf were friends who had collaborated extensively like in \cite{AlexandroffHopf}.
Textbooks cover it in \cite{JM,Spanier,GriffithsHarris,GD}. A version of the fixed point theorem 
in dimension $1$ is the
{\bf Nowakowski-Rival fixed edge theorem} \cite{NowakowskiRival} and a 
generalization \cite{Espinolaa} to a commutative family of such maps.
For the classical Lefschetz fixed point theorem and Brouwer fixed point
theorem see \cite{Lefschetz26a,BrownLefshetz}. The Brouwer fixed point
theorem is usually stated for balls in Euclidean space, where one can
prove it also with discrete methods, like Gale's proof using q-dimensional
Hex games. The text mentions some work on discrete Hurwitz \cite{TuckerKnill}
which I discussed in 2012/2013 with Thomas Tucker. I talked about this a couple of times, 
like at \cite{KnillBaltimore} or on July 7, 2026 in a youtube presentation. One reason
why not pursue this more was when realizing that the result is a rather direct
consequence of the Burnside lemma applied to each dimension sector of the 
simplicial complex. We hope to mention it more in a later volume of this project. \\
\index{Thomas Tucker}
\index{Nowakowski-Rival fixed point theorem}
\index{Brouwer fixed point theorem}

{\bf 8. Sphere formula} \\

The sphere formula is also a show case how to go from 
graphs to simplicial complexes. In both cases there are 
``unit spheres". For graphs, it is the graph generated by the neighbors. 
For simplicial complexes, the unit spheres are the boundaries of the 
smallest open sets. For traditional Euler characteristic, the 
sphere formula follows from the Green formula 
$\chi(G) = \sum_{x} \omega(x) \chi(U(x))$ which is ${\rm str}(L)$.
One also has $\chi(B(x))=1$ so that $\chi(G)=\sum_x \omega(x) \chi(B(x))$ 
Together with the definition $\chi(G) = \sum_x \omega(x)$ this proves the
sphere formula. The sphere formula immediately shows that
odd dimensional manifold have Euler characteristic zero. 
\cite{KnillEnergy2020,Sphereformula}. The Hydrogen relation was also discussed in
\cite{DehnSommerville} and \cite{AverageSimplexCardinality} where it is
shown that $g=\log(f)$ has the property that $g'(1)$ is the average 
simplex cardinality.  For Dehn-Sommerville, see
\cite{Klee1964,NovikSwartz,MuraiNovik,LuzonMoron,BrentiWelker,Hetyei,Klain2002,BergerLadder}.
Why unit spheres or hyperbolicity in simplicial complexes are less studied is maybe also 
due to notation. I myself got motivated from dynamical systems, especially classes of Lanford
or \cite{Lan85} or \cite{Yoc95}, the Hillerod conference I attended, as well as attemps of 
mine to understand hyperbolicity in ergodic set-ups \cite{Kni91,knill_1992}.
Traditionally, one has focused on the vertices of a simplicial complex
rather than all simplices. 
\cite{Veblen} defined a {\bf neighborhood} of a k-simplex, still using
geometric realizations.  
Veblen also used already the terminology of ``stars".
The notion of ``Stars" was established  in combinatorial topology like \cite{Alexandroff1937} .
It also entered some calculus textbooks like Whitney \cite{Whitney1957}. 
But even Alexandroff look at it primarily at geometric realizations of cell complexes.
For the unit sphere $S(v)$ of a zero dimensional complex, there
is no interesting hyperbolicity to be observed. There are strange names  also used in the
literature for $S(v)$. One calls it a ``link". Maybe, one should credit the discrete
Morse theory of Forman for centering on simplices as ``points". This is also the 
point of view for Barycentric refinement. If $G$ is a complex, then the 
Barycentric refinement of $G$ has the simplices of $G$ as ``points". The Barycentric 
refinement, then formalizes this and has the original simplices now zero dimensional objects. 
Most of the literature still looks at Barycentric refinements not combinatorially but by using
a geometric realization of the complex in an Euclidean space.  \\
\index{discrete Morse theory}
\index{link}
\index{Hydrogen relation}

{\bf 9. Level sets} \\

Much of calculus or algebraic geometry or differential topology deals with 
level sets of functions. Arthur Sard \cite{Sard1942,Sard1958} and Anthony Morse 
\cite{Morse1939} established that
a smooth map $f$ from a manifold $M$ to a manifold $N$ has the property that 
for almost all $y \in {\rm ran}(f)$, the level set $\{ x \in M, f(x) = y \}$ 
is a manifold. Already simple examples of polynomial maps like $x^2-y^2=0$ show
that level sets are in general not manifolds. It can come a bit as a shock
that in the discrete singularities do not occur. Level sets of any function
always are manifolds if not empty. Of course, we need to define level set properly.
This was done in \cite{KnillSard} in 2015. 
It was generalized in 2024 \cite{ManifoldsPartitions}. Working with Jennifer Gao
who wrote a thesis in that area helped a lot during this time \cite{Gao2024}. 
The level set result get for teaching purposes in Math 136
differential geometry classes. It is simple enough to be squeezed into a discussion
of what a manifold is. 
In \cite{DehnSommervilleManifolds}, it was generalized to Dehn-Sommerville 
manifolds, an approach to manifolds which is recursive but only makes an 
assumption about Euler characteristics. 
For the history of the development of the classical notions of manifolds, see
\cite{HistoryTopology,Scholz,Manolescu}.

\index{Jennifer Gao}

{\bf 10. Index formula} \\

The formula first appeared in \cite{indexformula} in 2012 and \cite{eveneuler} in 2013.
There was no level set construction yet and things were more complicated and ugly from
a current perspective. It was an application for \cite{indexexpectation} but it required
to define a level set structure defined by a function. It was valuable time however to 
work on this because it enabled the level set result to appear three years later. 
I like the index formula because it
is maybe the swiftest way to see why odd-dimensional manifolds are flat. This observation
was one of the earliest I made and it bothered me in 2010 that while possible to verify it
for 3-manifolds, I had been unable to verify it for 5 manifolds already. 
It is also intuitive to see curvature as a genus expectation of random sub-manifolds
in unit spheres. This shows that if we take a 4 manifold, a small sphere $S_r(x)$ in it and
a random function, then the expected genus of contour surfaces  $\{ f(y)=f(x)\}$ is
curvature. Since this is related to curvature of random 2-manifolds, this is
conceptionally closely related to scalar curvature which is the expected sectional 
curvature at a point. It motivates to look at Euler characteristic as a serious
functional with possible physical relevance. 
The index formula is non-trivial already for 2-manifolds $G$. Put random function values
on the vertices of a 2-manifold and at each node $v$ look at the average number $X(v)$
of sign changes of $f(w)-f(v)$ when $w$ goes along the circle around $v$. Then
$K(v) = 1-{\rm E}[X(v)]/2$ is equal to $1-d(v)/6$, the Eberhard curvature.  \\

{\bf 11. Quadratic cohomology} \\

This started in \cite{valuation} in 2016 and cohomology in \cite{CohomologyWu},
where the terminology "Wu cohomology" or "connection cohomology" was used.
Wen-Tsun Wu in 1953 first considered this multi-linear valuation in 1953 \cite{Wu1953}.
\cite{Gruenbaum1970} picked it up and conjectured the existence of Dehn-Sommerville
invariants, which was confirmed in \cite{valuation}.
Tamas Reti pointed out in 2016 a relation to the {\bf Zagreb index}.
When the fusion inequality for simplicial cohomology \cite{fusion1} was generalized
to higher cohomology, \cite{fusion2}, the name "quadratic cohomology" stuck.
We mentioned the theory last in \cite{ConnectionDirac} in the context of Birkhoff's
fixed point theorem.
I gave a table talk on March 8, 2016 about Wu characteristic.
\cite{KnillWuCharacteristic}.
It needs to be studied much more especially also in the context of the fusion 
story \cite{fusion1,fusion2},
where the cohomology is used also for open sets and interface cohomologies between
an open and closed set are used. One can look at the quadratic cohomology for the 
complement of a path in a 3-sphere for example. Does this give interesting 
knot invariants?  In \cite{CohomologyWu} we computed several examples. An interesting
one is the dunce hat because its quadratic cohomology is different from the quadratic
cohomology of a point. It demonstrates that quadratic cohomology is not a homotopy
invariant. 
\index{Wen-Tsun Wu}

{\bf 12. Higher Green}  \\

This section was covered \cite{CharacteristicTopologicalInvariants}, work done mostly 
in the winter break of 2022, when visiting Burlington in Vermont. It generalizes
the energy theorem in the case when $k=2$ and $m=1$, the Euler characteristic case. 
The case $k=1$ has appeared as ${\rm str}(g)=\chi(G)$. 
It is remarkable, how prominently Green functions appear in physics. One usually 
associates Green functions with the inverse of Laplacians but this is only the 2-particle
story. In the continuum, it leads to
technical difficulties as the Laplacians are not invertible. Particle interaction 
potentials have singularities.  The Green function of the Laplacian in Euclidean space 
$\mathbb{R}^3$ for example leads to the $1/r$ potential energy, we see in physics, like in 
electro magnetism or gravity.  Every geometric space defines so a potential, in $\mathbb{R}^2$ 
it would be $\log$ in $\mathbb{R}^4$, two celestial bodies would interact with a $1/r^2$ 
potential. The work on higher Green functions was motivated by physics. In particle physics,
one looks at {\bf $k$-point functions} $\langle \Omega,\phi(x_1) \dots \phi(x_k), \Omega \rangle$
or Green functions to calculate the out come of particle interactions. 
This is usually done using a series expansion and estimated using Feynman diagrams. But in 
the discrete, there are no infinities, everything is finite. The 1-point function $k=1$ case 
can be thought of as vacuum expectation. The 2-point functions  $k=2$ could have relations with
probability amplitudes for particles to go from $x_1$ to $x_2$. The general case $k$
could relate to scattering amplitudes of $k$ particles. It is important to note
however that all what is done here is just geometry and mathematics, actually combinatorics.
Physics is motivation however.

 \vfill \pagebreak
                         \section*{Structures}

\vspace{3mm}
A {\bf finite simple graph} is a pair $(V,E)$, where $V$ is a finite set and $E \subset \{ \{ a,b \}, a,b \in V, a \neq b \}$
is a set of pairs of points. Elements in $V$ are called {\bf vertices}. Elements in $E$ are called {\bf edges}. 
Given a subset $W$ of $V$, the {\bf graph induced by $W$} is the graph 
$(W,\{ (a,b) \in E, a \in W, b \in W \})$. For example, if $W$ is the set 
of vertices that are connected directly to a vertex $v$, it induces the {\bf unit sphere} $S(v)$. 
\index{finite simple graph}
\index{vertex}
\index{edge}
\index{unit sphere}
\index{induced graph}

\vspace{3mm}
A graph $\Gamma=(V,E)$ comes naturally with a simplicial complex $G$, where $G$ are the 
vertex sets of complete subgraphs of $\Gamma$. The {\bf maximal dimension} of $\Gamma$ is the 
maximal dimension of $G$. A triangle-free graph with at least one edge has dimension $1$. 
A graph without edges has dimension $0$. The empty graph has dimension $-1$. 
A tetrahedron $K_4$ has dimension $3$. The surface of 
the tetrahedron, is the 2-skeleton complex of $K_4$ consisting of all 
simplices of dimension $2$ or lower. If $G=\{ \{1\},\{2\},\{3\},\{1,2\},\{2,3\},\{3,1\},\{1,2,3\} \}$
is $K_3$, its 2-skeleton complex is 
$H=\{ \{1\},\{2\},\{3\},\{1,2\},\{2,3\},\{3,1\} \}$. A complex $G$ is called {\bf contractible}
if there exists $x$ such that $S(x)=\overline{U}(x) \setminus U(x)$ and $G \setminus U(x)$
are both contractible with the induction assumption that $1=\{ 1\}$ is contractible. 
The complex $K_5$ is a contractible complex of maximal dimension $4$, 
the $2$-skeleton complex $H$ is not contractible. 

\vspace{3mm}
A finite abstract simplicial complex $G$ defines a finite simple graph $\Gamma=(V,E)$ in which 
$V=G$ are the vertices and a pair $(a,b)$ is in the edge set $E$, if $a \subset b$ or $b \subset a$. 
The complex $G$ could also be defined to be contractible if its graph is contractible. 
For simplicial complexes, there are more contraction steps: $G$ is contractible, if there is 
a simplex $x \in G$ for which both the unit sphere $S(x)=\overline{U(x)} \setminus U(x)$ as well 
as $G \setminus U(x)$ are contractible simplicial complexes. If $G$ is the Whitney complex of a
graph $\Gamma=(V,E)$ and $x=v$ is a 0-dimensional simplex, then $G \setminus U(x)$ is the 
Whitney complex of $\Gamma \setminus v$.
\index{Graph from a complex}
\index{contractible complex}
\index{contractible graph}

\vspace{3mm}
If $(V,E)$ is a graph, we can look at the graph of its complex $G$. This is the
{\bf Barycentric refinement} of $(V,E)$. The Barycentric refinement of a complete
graph $K_{q+1}$ is a $q$-manifold with boundary. The Barycentric refinement of a 
triangle $K_3$ for example is a wheel graph with $7$ vertices. 
If $G$ is a complex, the Barycentric refinement is done by first looking at its 
graph $\Gamma$ and then taking the Whitney complex of that graph. 
In any case, the combination of producing a graph from a complex and producing a complex
from a graph is a Barycentric refinement. 
\index{Barycentric refinement}

\vspace{3mm}
Given two finite simple graphs $A,B$, the {\bf join} is the
graph $A \oplus B = (V,E)$, where $V=V(A) \cup V(B)$ and $E=E(A) \cup E(B) \cup \{ (a,b), a \in V(A), b \in V(b) \}$
are the edges. The simplicial complex $G$ of $A \oplus B$ is the join of the simplicial complexes
of $A$ and $B$. Every simplex $x \oplus y$ of $G$ is a join of two simplices $x$ and $y$. The $f$-vector of $G$ 
encodes the cardinalities $f_k=|G_k|$  of $k$ dimensional simplices in $G$. It defines the simplex generating 
function $f(t) = 1 + f_0 t + \cdots + f_q t^{q+1}$. We have $f_{A \oplus B}(t) = f_A(t) f_B(t)$. 
On the simplicial complex level, the join of two simplicial complexes $G,H$ is the union 
$G \cup H \cup \{ x \cup y, x \in G, y \in H \}$. If $G,H$ are the Whitney complexes of graphs, then their
join is the Whitney complex of the graph join. 
\index{join}

\vspace{3mm}
A {\bf finite set} is a collection $G$ of finite objects. It is a standard assumption
in set theory that a set does not record multiplicities of elements. The object $\{a,a,b,a,b\}$ 
for example is identified with $\{a,b\}$ if it is considered to be a set. 
If multiplicity should be needed, one could write it as a multi-set or a function 
$\{a,b\} \to \mathbb{N}$ encoding the multiplicities. If the order should matter, 
one could write it as a {\bf finite sequence} $\{ a,a,b,a,b \}$ or a {\bf word} $aabab$ over an 
alphabet.
\index{finite set}

\vspace{3mm}
If $V$ is a finite set, then $G=2^{V}$ is the {\bf power set} of $V$, the
set of all subsets of $V$. For $V=0$, one gets $G=2^0=1=\{0\}$, the 1-point set.
This set has a {\bf partial order structure}
$A \leq B$ in the form of $A \subset B$. It is also standard in set theory not to assume
the subset operation to be strict. The statement $A \subset A$ for example is always true. 
The set $2^V$ contains the empty set $\emptyset$. If $V=\{\} =\emptyset$ is the empty 
set, then $2^{V}=\{ \emptyset \} =1$ is a set with one element. 
Let $|A|$ denote the {\bf cardinality} of the set. The cardinality is a non-negative integer
$\mathbb{N}=\{0,1,2,3, \dots \}$. We would write $\mathbb{N}^+$ for the positive integers
and then there are the integers $\mathbb{Z}$, the group completion of $\mathbb{N}$. 
We have $|0|=0,|1|=1$ and $|2^A|=2^{|A|}$, justifying the notation. 
There is a natural identification of $2^{A+B}$ with the {\bf Cartesian product} $2^A \times 2^B$
because if $x$ is a subset of $A + B$, then $x$ defines $(a,b)$ with $a = x \cap A$ and 
$b=x \cap B$. But $(a,b) \in 2^A \times 2^B$. 
\index{power set}

\vspace{3mm}
There is a monoid operation on the category of finite sets. It is the {\bf disjoint union} $+$,
and is also known as {\bf co-product}. It is associative and so a {\bf monoid structure}.
The empty set $\emptyset$ is the zero element. If this monoid is group completed 
to become a group, one has then not only sets but ``negative sets".
If the co-product is {\bf group completed} with negative
sets, one has an additive group. 
\index{additive group of sets}
\index{monoid}
\index{co-product}

\vspace{3mm}
The {\bf Cartesian product} $A \times B$ is defined as the set of all ordered pairs $\{ (a,b), a \in A, b \in B\}$. 
The singleton set $1$ is the 1-element in this multiplication. The {\bf Cartesian ring} $(FinSet,+,0,\times,1)$ 
is a commutative ring with $1$ if one identifies $(-A) \times B = A \times (-B) = -(A \times B)$. 
The category FinSet
so naturally extends to a commutative ring in which $1$ is the one-element and where
$0$ is the zero-element. If the cardinality function is extended to negative sets as $|-A|=-|A|$, the
cardinality function defines a homomorphism from the ring of finite sets to $\mathbb{Z}$. If 
sets of the same cardinality are identified, then this is an isomorphism. 
\index{category of finite sets}

\vspace{3mm}
A {\bf finite abstract simplicial complex} $G$ is a finite set of non-empty sets that
is closed under the operation of taking non-empty subsets. 
This geometric structure has only one axiom, the {\bf hereditary axiom}. This simplicity
makes it appealing. There is a prise to pay however. Sets of sets are less intuitive than
graphs. 

\vspace{3mm}
Given a simplicial complex, one can assume that every $x \in G$ comes with an order. The choice 
of the order corresponds to picking a coordinate system and is irrelevant. It does not matter
what order we chose for each $x$, but it can be important to have an order in order to do 
linear algebra. We do not insist that the order is inherited by subsets. There is no order
on maximal simplices in general which works like that, the simplest example being the Moebius
strip complex, which is the simplicial complex of the graph complement of $C_7$. 
If $G$ contains subsets of integers, we can take for example the natural order. We also
assume that the elements of $G$ are ordered according to dimension to get block matrices. 
We for example write $G=\{ \{1\},\{2\}, \{1,2\}\}$ the complex of $K_2$ and 
not $G=\{ \{1,2\},\{1\},\{2\}$. Again, this is just a choice of labeling a basis when
doing linear algebra and irrelevant for any geometry.

\vspace{3mm}
As custom in set theory, sets contain different elements. 
$G=\{ \{1,2\},\{1\},\{2\} \}$ is an example of a simplicial 
complex. The structure $G=\{ \{1\}.,\{1,1\} \}$ would be identified with 
$\{ \{1\},\{1\} \}$ and so with $\{ \{1\} \}$.
A set does not record multiplicities. If this is needed, like for multi-graphs, we would
use sequences. For example, the sequence $( \{1\},\{2\},\{1,2\},\{1,2\} )$ is a quiver
with two edges connecting $1$ with $2$. 
\index{quiver}

\vspace{3mm}
An equivalent definition of simplicial complexes that is often
seen in the literature invokes the
{\bf vertex set} $V=\bigcup_{x \in G} x$. The assumption is then that
$G$ is a subset of $2^V \setminus \{ \{\} \}$ that does not contain the empty set
and that satisfies the heriditary axiom. In that frame work, a simplicial complex
is a pair $(V,G)$ where $G \subset 2^V$. 
\index{vertex set}
\index{hereditary axiom}

\vspace{3mm}
More general than a simplicial complex is a {\bf set of sets}. This structure has also been
called a {\bf hypergraph}. A hypergraph in general does not allow a reasonable calculus in 
the sense that one can not define a non-trivial exterior derivative.
A hypergraph also is allowed to contain the empty set. 
In order for the theory to work properly, it is important not to let the empty set to 
be in. The empty set can however be a structure itself. It is the zero element in all 
structures. The heriditary assumption is needed to get a notion of derivative that works. 
\index{hypergraph}

\vspace{3mm}
There are flavors of combinatorial topology that include $\emptyset$ into $G$.  
The {\bf Euler characteristic} $\chi(G) = \sum_{x \in G} \omega(x)$ with 
$\omega(x)=(-1)^{{\rm dim}(x)}$ for example would not work. 
The genus $1-\chi(G)$ is popular as it renders the join multiplicative.
\index{genus}
\index{Euler characteristic}

\vspace{3mm}
There are compelling reasons to keep $\chi(G)$ the Euler characteristic as we want compatibility
with addition (disjoint union) and multiplication (Cartesian product). There are also 
reasons to assume that simplices are contractible objects in the sense that its closure $\overline{\{x\}}$ 
is contracdtible. One should distinguish the element $x \in G$ with $K=\overline{\{x\}} \subset G$ even so they
can be naturally identified.  The object $K$ is simplicial complex of Euler characteristic $1$. 
And then there is $U=\{ x\}$ which is an open set in $K$. The complement $K \setminus U$ is the 
{\bf boundary complex} of $K$.
It is closed, a simplicial complex and always a sphere. 
Spheres are always non-contractible objects of Euler characteristic $1+(-1)^q$. 
The void $0$ is the $(-1)$ sphere. It is the boundary complex of the 1-point complex $1=\{ 1\}$.
It has Euler characteristic $0$ and is not contractible. 
\index{boundary complex}

\vspace{3mm}
Simplicial complexes can be defined from other structures. Any graph defines a {\bf Whitney complex}, in which 
the vertex sets of complete subgraphs are the sets. It is also called the {\bf flag complex} or {\bf order complex}.
The empty graph generates the void. It is not a complete complex as a simplicial complex never contains
the empty set. An other complex is the {\bf forest complex}, in which the edges of a forest are the sets. 
A forest decomposes into trees and trees are assumed in this set-up to be one-dimensional, meaning to 
have maximal dimension $1$. Single points = seeds, are not yet considered trees. In other contexts, like
for the matrix forest theorem, it is helpful to include seeds as trees. 
\index{Whitney complex}

\vspace{3mm}
Simplicial complexes can be derived from already given simplicial complexes.
The intersection of two simplicial complexes is a simplicial complex. If two simplicial complexes are
disjoint, then their intersection is the void, which is a simplicial complex. 
The union of two simplicial complexes $G=A \cup B$ is a simplicial complex as checking the definition shows.
If $x \in G$, then either $x \in A$ or $x \in B$ so that the hereditary axiom holds. 
If $G$ is a simplicial complex of dimension $q$, the {\bf skeleton complex} of dimension $k \leq q$ is the 
set of sets for which the dimension is smaller or equal than $k$. 
\index{skeleton complex}

\vspace{3mm}
Given an arbitrary set $A$ of non-empty sets, one can look at the {\bf closure}
$\overline{A}$ of $A$. It is the smallest simplicial complex containing $A$. 
To get $\overline{A}$, just adjoin all the non-empty subsets of $A$. 
The name closure is adequate because it is the closure operation in the Alexandroff
topology, the topology first considered by Pavel Alexandrov. 
\index{Alexandroff topology}

\vspace{3mm}
A {\bf topological space} $(X,\mathcal{O})$ is a set $X$
with a collection $\mathcal{O}$ of subsets of $X$ with the
property that both $\emptyset$ and $X$ are in $\mathcal{O}$ and that
$\mathcal{O}$ is closed under finite intersections and closed
under arbitrary unions. A topological space is {\bf finite}, if $X$
is a finite set. A {\bf base} $\mathcal{B}$ of a topological space is a subset
of $\mathcal{O}$ that generates $\mathcal{O}$ in the sense 
that every $A \in \mathcal{O}$ can be written as a union of
base elements, where the empty union $\bigcup_{A \in \emptyset} A=\emptyset$. 
The base does not have to contain the empty set $\emptyset$. 
\index{finite topology}
\index{topology}

\vspace{3mm}
The {\bf Alexandroff topology} is the topology generated
by the base $\{ U(x) \}_{x \in G}$, 
given by the {\bf stars} $U(x)=\{y \in G, x \subset y\}$. 
A {\bf sub-base} is $\{ U(v) \}_{x \in V=G_0}$ because intersections
of the sub-base produce base elements. If $G$ has positive maximal 
dimension, the topology is non-Hausdorff. Given an edge $(a,b)$, 
the points $a,b$ can not be separated by open sets as every open 
$U(a) \cap U(b)= U( (a,b) )$ is open and not-empty. 
\index{Alexandroff topology}

\vspace{3mm}
An arbitrary topology has the {\bf Alexandroff property}, if arbitrary 
intersections of open sets are open. Every finite
topology automatically has the Alexandroff property.
If the abstract simplicial complex is not finite 
but has a finite maximal dimension, then still, 
its topology has the Alexandroff topology. 
\index{Alexandroff property}

\vspace{3mm}
A topology $(G,\mathcal{O})$ is {\bf Kolmogorov} (T0) if for any two distinct points $x,y \in G$,
there exist neighborhoods $U(x),U(y)$ such that either $x \notin U(y)$ or 
$y \notin U(x)$. A topology is {\bf Fr\'echet} (T1) if for any two distinct points $x,y \in G$,
there exist neighborhoods $U(x),U(y)$ such that $x \notin U(y)$ and $y \notin U(x)$. 
A topology is {\bf Hausdorff} (T2) if for any two distinct points $x,y \in G$, 
there exist neighborhoods $U(x),U(y)$ such that $U(x) \cap U(y) = \emptyset$. 
\index{Kolmogorov property}
\index{Fr\'echet property}
\index{Hausdorff property}

\vspace{3mm}
The topology on a simplicial complex $G$ is always {\bf Kolmogorov}. Proof:
given $x \neq y \in G$. If there is an inclusion like $x \subset y$, then $U(y) \subset U(x)$ but $x \notin U(y)$. 
If there is no inclusion, then both $x \notin U(y)$ and $y \notin U(x)$. 
But a topology on $G$ is neither {\bf Fr\'echet}, nor {\bf Hausdorff} if its maximal
dimension is positive. If $x \subset y$, then $y$ is always in $U(x)$. 
As any finite topology, it is {\bf Alexandroff}, meaning that there are smallest
neighborhoods $U(x) = \{ y \in g, x \subset y\}$ of every $x \in G$. 
A finite topological space that is 
Hausdorff is necessarily the discrete topology, in which all subsets are open.
Finite topologies are never Hausdorff unless they are zero-dimensional. 

\vspace{3mm}
The Alexandroff topology of a simplicial complex is close to the {\bf Zariski topology} 
because the closed sets are the sub-simplicial complexes. We could 
associate sub-simplicial complexes of maximal dimension $q$ which have the property 
that all unit spheres have maximal dimension $q-1$ as some sort of variety.

\vspace{3mm}
Define the {\bf unit ball} $B(x) = \overline{U(x)}$ and
the {\bf unit sphere } $S(x) = B(x) \setminus U(x)$. The later is the boundary 
of the unit ball and is again closed as the intersection of a closed set $B(x)$ with 
a closed set $U(x)^c$. If $x=v$ is 0-dimensional, the
unit sphere $S(x)$ is sometimes called the {\bf link}. 

\vspace{3mm}
Note that $G$ is a set of sets and that elements in $\mathcal{O}$ are a set of subsets 
of $G$. The open sets in $G=\{ \{1\},\{2\},\{1,2\} \}$ 
are $\{ \emptyset, \{ \{1\},\{1,2\} \}, \{ \{2\},\{1,2\} \}, \{ \{1,2\} \}, G \}$.
So, it is important to distinguish a simplex $x \in G$ with $\{ x \}$ which is 
in general neither open nor closed. Maximal simplices $x$ produce the open set $U(x) = \{ x\}$, 
vertices $x$ produce closed sets $C(x) = \{x\}$. If $x$ is an intermediate dimension, 
then $\{x\}$ is neither open nor closed. The smallest closed set containing $x$ is 
$\overline{ \{ x \}}$, the smallest open set containing $x$ is $U(x)$. 

\vspace{3mm}
One could also look at topologies on the vertex set $V$ but that would be different. 
For example, if $G$ is the complex of $K_4$ with vertex set $V=\{1,2,3,4\}$, then
the pseudo circle $\mathcal{O} = \{ (1,2,3,4),(1,2,3),(1,2,4),(1,2),(1),(2),(0) \}$ 
a subset of $G$, not a set of sets from $G$. This example is known as the {\bf small circle}.
It is neither an open set nor a closed set in $G$.

\vspace{3mm}
If $G$ is a finite abstract simplicial complex, 
a {\bf contraction step} is defined as a removal of a star $U(x)$
with contractible boundary $S(x)=\delta U(x)$ from $G$. Formally, this means
$G \to G \setminus U(x)$ under the assumption that 
$S(x) = \overline{U(x)} \setminus U(x)$ is contractible. 
The inverse operation is called an {\bf expansion step}. Two complexes 
are called {\bf homotopic} if there exists a finite number of homotopy expansion or 
contraction steps, deforming one to an other. 
\index{homotopy contraction}
\index{contraction}
\index{expansion}

\vspace{3mm}
If $G$ can be reduced to $1$ by a finite number of contraction 
steps, the complex $G$ is called {\bf contractible}. This is 
very different from being  homotopic to $1$. The {\bf dunce hat} is an 
example of a complex that is homotopic to $1$ but not contractible. 

\vspace{3mm}
If $G$ is contractible and $H$ is arbitrary, then $G \oplus H$ is 
contractible. For example $G \oplus 1$, the cone extension is always
contractible. Every ball $B(x) = \overline{U(x)}$ is contractible. 
Contractible simplicial complexes form a {\bf semi-group} within the 
monoid of all simplicial complexes. It is not a sub-monoid because the 
zero element is not considered contractible. 

\vspace{3mm}
No q-manifold is contractible unless it is the zero-dimensional complex
$G=1$, the one point complex. No sphere is contractible: the empty set which is 
the $(-1)$-sphere is not contractible by definition. Recursively, no
$q$-sphere for $q>0$ is contractible.  The semi-group of contractible 
complexes is disjoint from the sphere monoid. 

\vspace{3mm}
Contractibility can be computed in a time complexity that grows
maximally quadratically in the number of elements of $G$: just go 
through all the $x \in G$ and check in each case whether $S(x)$ and
$G \setminus U(x)$ are both contractible. Both 
simplicial complexes $S(x)$ and $G \setminus U(x)$ are smaller. 
\index{complexity of contractibility}

\vspace{3mm}
On the other hand, there is no fixed Turing machine that takes as 
an input a simplicial complex $G$ and tells whether $G$ is homotopic
to $1$ or not. The distinction between contractibility and {\bf homotopic to 1}
could not be bigger. This is a reason not to use words like {\bf collapsibility}
and {\bf contractibility} to make the distinction since both terms are in 
colloquial language associated with a shrinking process. But homotopic to 1
can mean that we have to take a complex and blow it up again and again until 
it eventually can be shrunk to a point.

\vspace{3mm}
Fix a finite abstract simplicial complex $G$. 
A linear functional $X$ on the set of $f$-vectors $\phi(G) = X \cdot f(G)$
is called a {\bf valuation} on $G$.
It satisfies $\phi(A \cup B) + \phi(A \cap B) = \phi(A) + \phi(B)$ for 
closed subsets. 
\index{Valuation}

\vspace{3mm}
The Hadwiger theorem statement that the set of 
valuations on a complex of maximal dimension $q$ has
dimension $q+1$ follows from the definition. A basis is 
given by the counting valuations $f_k(G)$. An example of a valuation is 
the Euler characteristic $\chi(G) = f_0-f_1+f_2 \cdots (-1)^q f_q$. 

\vspace{3mm}
A valuation $X$ that satisfies $X(A)=1$ for every complete sub-complex $A$ of $G$
is $\chi(G)$.  Proof: If $G$ has maximal dimension $q$, then the space of valuations has dimension $q+1$.
Valuations are of the form $X(A) = X \cdot f(A)$.
If $X(A)=1$ for every complete graph, then $X \cdot f(K_k) = 1$, for the
$f$-vectors $f(K_k)$ of $K_k$. But since these vectors form a basis in the
vector space $\mathbb{R}^{d+1}$, we have uniqueness and $X=\chi$.

\vspace{3mm}
Let $(G,G)$ be the set of all pairs of simplices $(x,y)$ with 
$x,y \in G$ such that $x \cap y$ is non-empty. 
More generally, let $(A,B)$ be the set of all pairs of simplices
$(x,y)$ with $x \cap y \neq 0$ and $x \in A, y \in B$. 
The {\bf f-matrix} of $(A,B)$ is $f_{k,l}(A,B)$, the number of elements
$(x,y)$ in $(A,B)$ with $x$ of dimension $k$ and $y$ of dimension $l$. 
\index{f-matrix}

\vspace{3mm}
The f-matrix defines a multi-linear valuation. The vector space of quadratic
valuations has dimension $(q+1)(q+2)/2$ as the numbers $f_{k,l}$ are symmetric,
this is the dimension of the space of symmetric $(q+1) \times (q+1)$ matrices. 

\vspace{3mm}
Valuations generalize counting or measure with some sort of symmetry. 
In the continuum, one restricts the theory of valuations to finite unions of
convex sets in $\mathbb{R}^n$.   
The theory of measures does not look at the internal structure 
of the sets which are measured, unlike the theory of valuations which look inside. 
In the continuum, one needs to look at the internal k-dimensional structure of a set using 
probabilistic tomography, like the Crofton formula. One can for example look at the length
of a convex set in $\mathbb{R}^2$ as an expectation of length of intersections with random 
lines. This satisfies the valuation property because for every experiment, we have the
valuation property. 

\vspace{3mm}
A reasonable {\bf geometry on a finite set} need a {\bf category} that enables
basic operations from calculus. We especially like to build products, quotients
or level sets. Any of these operations extends the scope of simplicial complexes
in general making it necessary to work in a larger category.

\vspace{3mm}
The category $Set$ of finite sets is an {\bf elementary topos}:
cartesian product,inclusion, exponentials and subobject classification works.
There is a natural generalization of this in the form of {\bf delta sets}.
\index{elementary topos}

\vspace{3mm}
A finite category $C$ defines the pre-sheaf category $Pre(C)=Set^{C^{op}}$.
{\bf Delta sets} are the pre-sheaf category of the simplex category $C$ with strict
inclusions as morphisms.
\index{delta set}

\vspace{3mm}
A delta set is a sequence of sets $G_0,G_1, \dots, G_q$ with {\bf face maps} $d_i^k: G_{k+1} \to G_k, i=0, k+1$
that satisfy $d_i^k d_j^{k+1}=d_{j-1}^k d_i^{k+1}$ for $i<j$. A simplicial complex is a delta set
The simplest example of a delta set that is not a simplicial complex is $G_0=\{0\},G_1=\{ (0,0) \}$
with maps $d_0( (0,0) ) = 0, d_1( (0,0) ) = 0$. Open sets in a simplicial complex are delta sets,
level sets of functions are delta sets. Products of simplicial complexes are delta sets.

\vspace{3mm}
An {\bf abstract delta set} is a finite set $G$ of n sets with a dimension function $R:G \to \{0,\dots, q\}$
and a matrix $D=d+d^*$, where $df$ is supported on $G_{k+1}=\{ R=k+1 \}$ if $f$ is supported on $G_k$.
Any delta set defines an abstract delta set. Abstract delta sets are sets of non-empty sets
with a notion of dimension and derivative. While one can avoid delta sets by going to a refinement and 
look at the Whitney complex of the graph defined by the delta set, there are computational advantages
to leave the delta set structure. 

\vspace{3mm}
The Dirac operator $D$ of a Cartesian product for example is 
much smaller than the Dirac operator of the Barycentric refined case. Similarly, the cohomology of a
level set in a simplicial complex can be determined faster than using the Whitney complex of the graph 
of this structure. A third example is looking at quotients $G/A$ where $A$ is a subgroup of the automorphism
group of a simplicial complex $G$. The quotient is often no more a simplicial complex.
\index{abstract delta set}

 \vfill \pagebreak

\printindex

\bibliographystyle{plain}

\end{document}